\documentclass[12pt]{amsart}

\usepackage{graphicx, epsfig}
\usepackage{amsfonts,amscd, amssymb,  epsf, epsfig, graphicx}
\usepackage{amssymb,amsmath,amsfonts,array,delarray,amsthm,graphics,graphicx}

\newcommand{\abs}[1]{| #1 |}

\def\bcases{\begin{cases}}

\def\ecases{\end{cases}}

\newcommand{\norm} [1]{\left\| #1\right\|}

\newcommand{\bea}{\begin{eqnarray*}}
\newcommand{\eea}{\end{eqnarray*}}

\newtheorem{thm}{Theorem}

\theoremstyle{definition}

\theoremstyle{remark}

\newcommand{\be}{\begin{equation}}
\newcommand{\ee}{\end{equation}}

\title{Properties of solutions of differential equations on foliations}
 
\begin{document}
\maketitle

\centerline{Xiaoai Chai\footnote{The first author is supported by REU grant at University of Michigan.
2011 AMS classification. Primary: 57R30; Secondary: 32S65} \and John E. Forn\ae ss\footnote{The second author is supported by an NSF grant DMS-1006294. }}

\begin{abstract}
In this paper, we discuss foliations by real curves. We investigate differential equations which are modifications of $\frac{du}{dx}=v$ along leaves. Our focus is on having a solution operator so that $u$ is continuous (or $C^\infty$) if $v$ is continuous (or $C^\infty$).
\end{abstract}

\section{Introduction}

In this paper, we study properties of solutions of differential equations that are defined on various real foliations. 

In section 2, we investigate $C^\infty$ foliations on compact sets $K$ defined by vector fields $X$. We study the differential equation $U+XU=V$ which is defined on the foliation. We show:
\begin{thm}
If $v$ is $C^\infty$ on a neighborhood of $K$, then there is a $C^\infty$ function $u$ on $K$ so that $U+XU=V$
\end{thm}
We prove Theorem 1 by mathematical induction on the number of derivatives. We also construct some particular examples to examine this theorem and to motivate the choice of equation.

In section 3, we study foliations with singularities. On the real line, we find a solvable differential equation and we prove that $u(x)$ is continuous and bounded if $v(x)$ is continuous and bounded. However, on a circle, we show that there is no such solvable differential equation.

In section 4, we construct certain line bundles on $C^\infty$ foliations. We show that in this case we can solve $\frac{du}{ds}=v$ with values in the line bundle. This is the main result in this paper:

\begin{thm}
If $v$ is $C^\infty$ on a foliated set $K$, then there is a bounded and continuous function $u$ such that $du/dt=v$.
\end{thm}

\section{$C^\infty$ Foliations}

\subsection{Real Line}
Consider a function $v(x)$ on a real line which is continuous and bounded. We want to find a bounded function $u(x)$ so that $v(x)=\frac{du(x)}{dx}$. If we let $v(x)=1$, then $u(x)=x+C$ where $C$ is a constant. Obviously $u(x)$ is not a bounded function. To solve this problem, we modify the differential equation by adding $u(x)$ to the left part of the equation and get $u(x)+\frac{du}{dx}=v(x)$. Then, we solve this new equation. First, we multiply the equation by $e^x$ and get $e^xu(x)+e^x\frac{du}{dx}=e^xv(x)$ which is the same as $\frac{d(e^xu(x))}{dx}=e^xv(x)$ and hence $\int  e^xv(x)dx=e^xu(x)$. We get a general formula for $u(x)$:
\begin{equation*}
u(x)=e^{-x} \int_{-\infty}^x v(t)dt.
\end{equation*}

Now we will prove the continuity of the function $u(x)$ by showing that $e^xu(x)$ is continuous if $v(x)$ is continuous.

Let $\epsilon$ be an arbitrary small positive number and $x\in \mathbb R$. We will show that there exist a real number $\delta$ such that $|e^xu(x)-e^{x'}u(x')|<\epsilon$ holds for every pair of $(x, x')$ if $|x-x'|<\delta$.
\begin{eqnarray*}
&&|e^xu(x)-e^{x'}u(x')|\\
&=& |\int_{-\infty}^x v(t)dt-\int_{-\infty}^{x'} v(t)dt|\\
&=& |\int_{x'}^{x} v(t)dt|
\end{eqnarray*}
Since $v(t)$ is bounded and continuous, there exist a $\delta$ such that $|e^xu(x)-e^{x'}u(x')|<\epsilon$ if if $|x-x'|<\delta$. This proves the continuity of $u(x)$.

\subsection{Torus}
\begin{figure}
\begin{center}
	\includegraphics[width=1.00\textwidth]{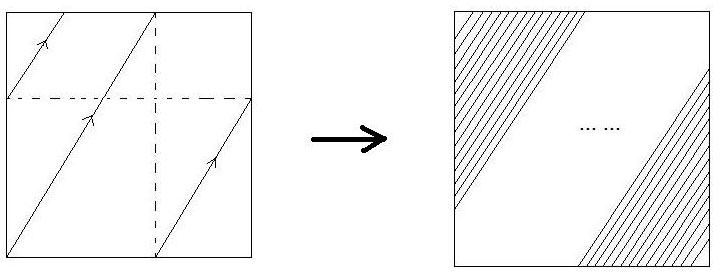}
	\label{Figure 1.}
\end{center}
\end{figure}

Now we expand this problem to a square on a plane with lines with an irrational slope $\sqrt{2}$. As shown in Figure 1 (a), whenever the line hits the upper (or right) boundary of the square, it goes straight down (or left) to the same x (or y) value on the lower (or left)  boundary with the same slope. Then we have a continuous and bounded function $v(x,y)$ and a bounded function $u(x,y)$ which is differentiable in the direction$(1, \sqrt{2}) $such that $v(x,y)=D_{(1,\sqrt{2})}u(x,y)$. We found that the same problem exists here as in the real line case condition. Let $x=t, y=\sqrt{2}t+C$, then we have $v(x,y)=v(t,\sqrt{2}t+C)$ and $ v(t,\sqrt{2}t+C)=\frac{du(t,\sqrt{2}t+C)}{dt} $. Setting$ v(x,y)=1$, we get $u(t,\sqrt{2}t+C)=t+C$ which is an unbounded function. Again we multiply both sides of the equation by $e^x$ and get $e^xv(t,\sqrt{2}t+C)=\frac{d(e^xu(t,\sqrt{2}t+C))}{dt}$ and $\int  e^xv(t,\sqrt{2}t+C)dt=e^xu(t,\sqrt{2}t+C)$. Here if we set $v(t,\sqrt{2}t+C)=1$, we get $u(t,\sqrt{2}t+C)=1$ which is a bounded function. 

However, we notice that by multiplying $e^x$ to the equation, we created a difference between the y values of the line when it reaches the right side of the square and moves to the left side. At the right side of the square, the line has the expression of $e^1v(x)$ while at the left side,  the expression becomes $e^0v(x)$. The y values of the two ends are not equal to each other any more. 
We use a line bundle to fix this difference. (We will discuss line bundle in more detail in section 4.) We connect the upper and lower sides together to form a hollow cylinder. Then, we bend the cylinder along its axis and make it to be a ring as shown in Figure 2. 
\begin{figure}[h]
	\centering
		\includegraphics[width=0.50\textwidth]{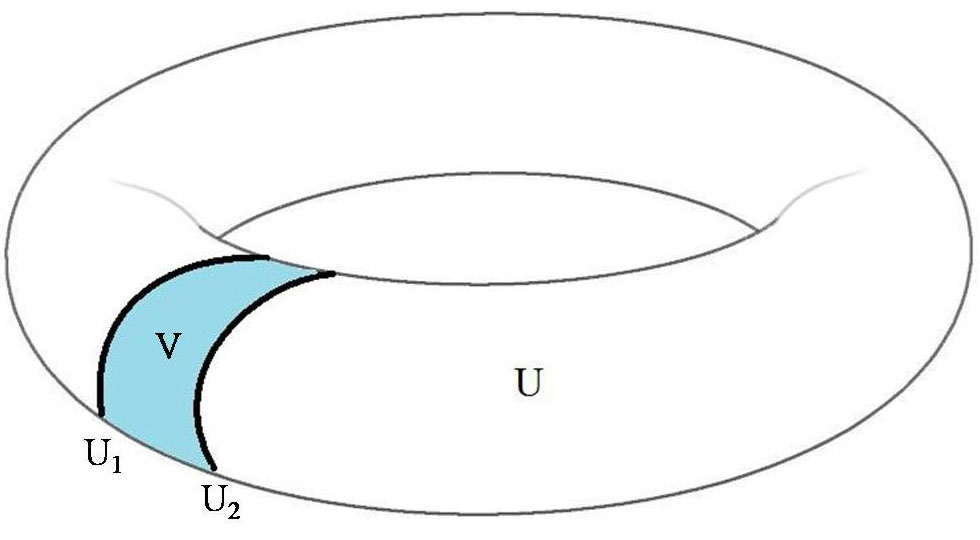}
	\label{fig:figure2}
\end{figure}

Let $n$ be a real number closed to 0 and $0 \leq n \leq 1$, U be a set that contains all the points on the square with $n \leq x \leq 1$, and V be a set contains the points with $0 \leq x \leq n$. We define $U_1$ as a set of points lies on $x=n$ and $U_2$ as the set of points lies on $x=0$. Both $U_1$ and $U_ 2$ are intersection of sets $U$ and $V$.
Based on the assumptions, we have       
\begin{equation}
U \times \mathbb{R} = \{(a, r)\}, a \mbox{ is a point and } a \in U \mbox{, } r \in \mathbb{R}
\end{equation}
\begin{equation}
V \times \mathbb{R} = \{(b, s)\}, b \mbox{ is a point and } b \in V \mbox{, } s \in \mathbb{R} 
\end{equation}
\begin{equation}    
U_1 \times \mathbb{R} = \{(a, r)| a=b, r=s \},
\end{equation}
\begin{equation}
U_2 \times \mathbb{R} = \{(a, r)| a=b, r=se \}.
\end{equation}
With $r=e^x$, we have $ \{(a, e^x)\} \in U \times \mathbb{R}$. For $U_2$, $ \{(a, e^{x-1})\} \in U_2 \times \mathbb{R}$. \\
Now, we are going to find the general form of $u(x)$. 
\begin{align*}
\frac{d}{dt} (e^tu(t,\sqrt{2}t+C))=e^tv(t,\sqrt{2}t+C),\\
e^tu(t,\sqrt{2}t+C)=\int_{-\infty}^t e^hv(h,\sqrt{2}h+C)dh,\\
u(t,\sqrt{2}t+C)=e^{-t}\int_{-\infty}^t e^hv(h,\sqrt{2}h+C)dh.
\end{align*}
Fix a point on the line with coordinates $(x_0, y_0)$, it is easy to know that $y_0=\sqrt{2}x_0+C$. According to the formulae we derived above, we have
\begin{equation}
u(x_0, y_0)=e^{-x_0} \int_{-\infty}^{x_0} e^hv(h,\sqrt{2}h+y_0-\sqrt{2}x_0)dh
\end{equation}
Since the function $v(x,y)$ is defined on the interval [0,1], in order to do the integration from negative infinity to $x_0$, we need to expand its definition. We assume that $v(x,y)$ is a periodic function with the period of 1,i.e. $v(x_0,y_0)=v(x_0+1,y_0)$. What we want to know is that whether the function $u(x,y)$ has the same property or not. \\
>From equation (5),
\begin{equation*}
u(x_0+1,y_0)=  e^{-(x_0+1)}\int_{-\infty}^{x_0+1} e^hv(h,\sqrt{2}h+y_0-\sqrt{2}(x_0+1))dh
\end{equation*}
Let $h'=h-1$,
\begin{eqnarray*}
u(x_0+1,y_0) &=& e^{-(x_0+1)}\int_{-\infty}^{x_0} e^{h'+1}v(h'+1,\sqrt{2}(h'+1)+y_0-\sqrt{2}(x_0+1))dh' \\
                         &=&  e^{-(x_0+1)}\int_{-\infty}^{x_0} e^{h'}ev(h'+1,\sqrt{2}h'+y_0-\sqrt{2}x_0)dh' \\
                         &=&  e^{-x_0}\int_{-\infty}^{x_0} e^{h'}v(h',\sqrt{2}h'+y_0-\sqrt{2}x_0)dh' \\
                         &=& u(x_0,y_0).
\end{eqnarray*}
So $u(x,y)$ is a periodic function with period of 1.
\\Next, we will prove the continuity of function $e^xu(x,y)$.

Let $\epsilon$ be an arbitrary small number. We will show that there exist a real number $\delta$ such that $v(h,\sqrt{2}h+y-\sqrt{2}x)-v(h,\sqrt{2}h+y'-\sqrt{2}x')<\epsilon/(2e^x)$ and $\abs{e^xu(x,y)-e^xu(x',y')} < \epsilon$ holds for every pair of $(x, x')$ and $(y, y')$ if $\abs{x-x'} < \delta$ and $ \abs{y-y'} < \delta$.
Because of the continuity of the function $v(x, y)$, there exists a $\delta_1$ such that $v(h,\sqrt{2}h+y-\sqrt{2}x)-v(h,\sqrt{2}h+y'-\sqrt{2}x')<\epsilon/(2e^x)$.\\
Without loss of generality, we assume $x>x'$.
\begin{eqnarray*}
& & \abs{e^xu(x,y)-e^xu(x',y')}\\
& = & \abs{\int_{-\infty}^xe^hv(h,\sqrt{2}h+y-\sqrt{2}x)dh-\int_{-\infty}^{x'}e^hv(h,\sqrt{2}h+y'-\sqrt{2}x')dh}\\
& = & \abs{\int_{-\infty}^xe^h[v(h,\sqrt{2}h+y-\sqrt{2}x)-v(h,\sqrt{2}h+y'-\sqrt{2}x')]dh\\
& &+ \int_x^{x'}e^hv(\sqrt{2}h+y'-\sqrt{2}x')dh}\\
& < &\abs{\int_{-\infty}^xe^h \frac{\epsilon}{2e^x} dh+\int_x^{x'}e^hv(h, \sqrt{2}h+y'-\sqrt{2}x')dh}\\
& = & \abs{ \frac{\epsilon}{2e^x} e^x +\int_x^{x'}e^hv(h, \sqrt{2}h+y'-\sqrt{2}x')dh}\\
& \leq & \frac{\epsilon}{2} + \abs{\int_x^{x'}e^hv(h, \sqrt{2}h+y'-\sqrt{2}x')dh}
\end{eqnarray*}
To make the condition $\abs{e^xu(x,y)-e^xu(x',y')} < \epsilon$ holds, we have $\abs{\int_x^{x'}e^hv(h, \sqrt{2}h+y'-\sqrt{2}x')dh} \leq \epsilon /2$. Let $k=e^hv(h, \sqrt{2}h+y'-\sqrt{2}x')$. Again because of the continuity of the function $v(x,y)$, there exist a $\delta_2$ such that 
\begin{eqnarray*}
\abs{\int_x^{x'}e^hv(h, \sqrt{2}h+y'-\sqrt{2}x')dh} & = &  (x-x')k\\
& = & \delta_2k \\
& < &  \frac{\epsilon}{2}
\end{eqnarray*}
\begin{equation}
\delta_2 < \frac{\epsilon}{2k}
\end{equation}
Let $\delta=\min \{ \delta_1, \delta_2 \} $. If $\abs{x-x'} < \delta$ and $ \abs{y-y'} < \delta$, then all conditions hold and this proves the continuity of the function $u(x, y)$.\\

\subsection{Two circles and a spiral between them}
Consider two concentric circles with radius 1 and 2 respectively. Between them, there is a spiral which approximates both of the circles but never touches them (Shown as Figure 3). 
\begin{figure}[h]
	\centering
		\includegraphics[width=0.40\textwidth]{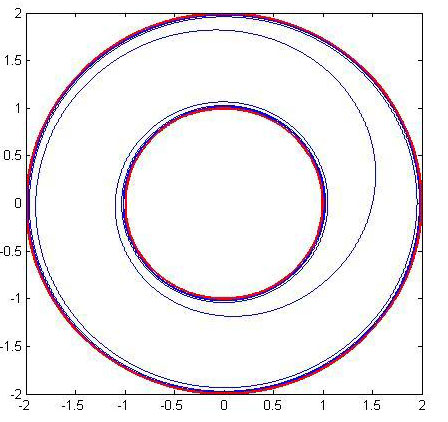}
	\label{fig:figure3}
\end{figure}

A continuous and bounded function $v(r, \theta)$ is defined on the two circles and the spiral.
 We want to find another continuous, bounded and differentiable function $u(r, \theta )$ such that $v(r, \theta) = du(r, \theta)/d\theta$. \\
On circles, we have
\begin{equation}
v(1, \theta)=v(1, \theta+2\pi),
\end{equation}
\begin{equation}
v(2, \theta)=v(2, \theta+2\pi),
\end{equation}
\begin{equation}
u(1, \theta) = u(1, \theta +2\pi), 
\end{equation}
\begin{equation}
u(2, \theta) = u(2, \theta +2\pi).
\end{equation}
On the spiral, the radius is a function of the angel which is 
\begin{equation}
r(\theta) =\frac{3}{2} + \frac{1}{\pi} \arctan \theta. 
\end{equation}
The function $v(r(\theta), \theta)$ has following property on the spiral:
\begin{equation}
\lim_{\theta \rightarrow -\infty} \abs{ v(r(\theta), \theta)-v(1, \theta)} = 0,
\end{equation}
\begin{equation}
\lim_{\theta \rightarrow +\infty} \abs{ v(r(\theta), \theta)-v(2, \theta)} = 0.
\end{equation}
We want to find that function $u(r(\theta), \theta)$ has the same property:
\begin{equation}
\lim_{\theta \rightarrow -\infty} \abs{ u(r(\theta), \theta)-u(1, \theta)} = 0,
\end{equation}
\begin{equation}
\lim_{\theta \rightarrow +\infty} \abs{ u(r(\theta), \theta)-u(2, \theta)} = 0.
\end{equation}
If we set $v(r(\theta), \theta)=1$, then $du(r(\theta), \theta)/d\theta=1 $ and hence $u(r(\theta), \theta) = \theta +C$ which is not a bounded function. This indicates that the differential equation $du(r(\theta), \theta)/d\theta=v(r(\theta), \theta)$ is not solvable. To solve the problem, we add $u(r(\theta), \theta)$ to the left side of the equation and get a new equation $u(r(\theta), \theta) + du(r(\theta), \theta)/d\theta=v(r(\theta), \theta)$. We found that we can solve this new equation by multiplying $e^\theta$ to both sides of the equation and get,
\begin{equation*}
e^\theta u(r(\theta), \theta) + e^\theta \frac{d}{d\theta}u(r(\theta), \theta)=e^\theta v(r(\theta), \theta),
\end{equation*} 
which is the same as 
\begin{equation*}
 \frac{d}{d\theta}e^\theta u(r(\theta), \theta))=e^\theta v(r(\theta), \theta).
 \end{equation*}
Now we are going to find the general form for $u(r(\theta), \theta)$.\\
On the circles, we have
\begin{equation*}
\frac{d}{d\theta}e^{\theta} u(r, \theta) \mbox{,  r=1,2.}
\end{equation*}
\begin{equation*}
u(r, \theta)=e^{-\theta} \int_{-\infty}^{\theta} e^hv(r, h)dh \mbox{, r=1,2.}
\end{equation*}
Fix a point on a circle with coordinates $(r_0, \theta_0)$, according to the formulae we derived above, we have
\begin{equation}
u(r_0, \theta_0)=e^{-\theta_0} \int_{-\infty}^{\theta_0} e^hv(r_0, h)dh \mbox{, r=1,2.}
\end{equation}
Then, we need to make sure that the function $u(r, \theta) $ is periodic with period of $2\pi$, i.e. we need to verify that the equation $u(r_0, \theta_0+2\pi) = u(r_0, \theta_0)$ holds.\\
From equation (16),
\begin{equation*}
u(r_0, \theta_0+2\pi)=e^{-(\theta_0+2\pi)} \int_{-\infty}^{\theta_0+2\pi} e^hv(r_0, h)dh
\end{equation*}
Let $h' = h - 2\pi$,
\begin{eqnarray*}
u(r_0, \theta_0+2\pi) & = & e^{-\theta_0+2\pi} \int_{-\infty}^{\theta_0} e^{h'+2\pi}v(r_0, h'+2\pi)dh\\
                                      & = & e^{-\theta_0} \int_{-\infty}^{\theta_0} e^{h'}v(r_0, h')dh\\
                                      & = & u(r_0, \theta_0)
\end{eqnarray*}
This proves that $u(r, \theta)$ is a periodic function with period of $2\pi$.\\
Let  
\begin{equation}
U(\theta) = u(r(\theta), \theta)
\end{equation}
\begin{equation}
V(\theta) = v(r(\theta), \theta)
\end{equation}
On the spiral, we have,
\begin{equation*}
\frac{d}{d\theta} e^{\theta} U(\theta)=e^{\theta} V(\theta),
\end{equation*}
\begin{equation*}
U(\theta)=e^{-\theta} \int_{-\infty}^{\theta} e^hV(h)dh.
\end{equation*}
Substitute (17), (18) back to the solution, we get
\begin{equation*}
u(r(\theta), \theta) = e^{-\theta} \int_{-\infty}^{\theta} e^hv(r(h), h)dh.
\end{equation*}
Fix a point on the spiral with coordinates $(r(\theta_0), \theta_0)$, according to the formulae we derived above, we have
\begin{equation}
u(r(\theta_0), \theta_0) = e^{-\theta_0} \int_{-\infty}^{\theta_0} e^hv(r(h), h)dh.
\end{equation}
Now we show the continuity of the solution we got. We are going to examine the continuity of $e^\theta u(r_0, \theta)(r_0=1,2)$, $e^\theta u(r(\theta), \theta)$, and the continuity between them.

We will fist show that the solution on the circles is continuous.
Let $\epsilon$ be an arbitrary small number. We will show that there exist a real number $\delta_1$ such that $\abs{e^\theta u(r_0, \theta)-e^{\theta'} u(r_0, \theta')} < \epsilon$ holds for every pair of $(\theta, \theta')$ if $\abs{\theta-\theta'} < \delta_1$.
Without loss of generality, we assume $\theta>\theta'$.
\begin{eqnarray*}
& &\abs{e^\theta u(r_0, \theta)-e^{\theta'} u(r_0, \theta')} ,     (r_0=1,2)\\
& = & \abs{ \int_{-\infty}^{\theta} e^hv(r_0, h)dh-\int_{-\infty}^{\theta'} e^hv(r_0, h)dh}\\
& = & \abs{\int_{\theta'}^{\theta} e^hv(r_0, h)dh}
\end{eqnarray*}
Because of the continuity of $v(r_0, h)$, there exists a $\delta_1$ such that 
\begin{equation*}
\abs{e^\theta u(r_0, \theta)-e^{\theta'} u(r_0, \theta')} < \epsilon
\end{equation*}
Then we show that the solution on the spiral is continuous.
Let $\epsilon$ be an arbitrary small number. We will show that there exist a real number $\delta_2$ such that $\abs{e^\theta u(r(\theta), \theta)-e^{\theta'} u(r(\theta'), \theta')} < \epsilon$ holds for every pair of $(\theta, \theta')$ if $\abs{\theta-\theta'} < \delta_2$.
Without loss of generality, we assume $\theta>\theta'$.
\begin{eqnarray*}
& &\abs{e^\theta u(r(\theta), \theta)-e^{\theta'} u(r(\theta'), \theta')} \\
& = & \abs{ \int_{-\infty}^{\theta} e^hv(r(h), h)dh-\int_{-\infty}^{\theta'} e^hv(r(h), h)dh}\\
& = & \abs{\int_{\theta'}^{\theta} e^hv(r(h), h)dh}
\end{eqnarray*}
Because of the continuity of $v(r(h), h)$, there exists a $\delta_2$ such that 
\begin{equation*}
\abs{e^\theta u(r(\theta), \theta)-e^{\theta'} u(r(\theta'), \theta')} < \epsilon
\end{equation*}
Then we show that the solution between the circles and the spiral is continuous.
Let $\epsilon$ be an arbitrary small number. We will show that there exist a real number $\delta_3$ such that $\abs{[u(r(\theta), \theta)- u(1, \theta')]} < \epsilon$ and $\abs{u(r(\theta), \theta)- u(2, \theta')}< \epsilon$ holds for every pair of $(\theta, \theta')$ and $(r(\theta), r_0)$ if $\abs{\theta-\theta'} < \delta_3$ and $\abs{r(\theta) -r_0}<\delta_3$.$( r_0 =1,2)$

Firstly, we show the continuity between $u(r(\theta), \theta)$ and $u(1, \theta)$ when $\theta \rightarrow -\infty$.
Without loss of generality, we assume $\theta > \theta'$.
\begin{eqnarray*}
& &\abs{ u(r(\theta), \theta)-u(1, \theta')}  \\
& = & \abs{e^{-\theta} \int_{-\infty}^{\theta} e^hv(r(\theta), h)dh-e^{-\theta'}\int_{-\infty}^{\theta'} e^hv(1, h)dh} \\
& = &  \abs{e^{-\theta} \int_{-\infty}^{\theta} e^hv(r(\theta), h)dh-e^{-\theta'}\int_{-\infty}^{\theta'} e^hv(1, h)dh\\
& &+ e^{-\theta} \int_{-\infty}^{\theta} e^hv(1, h)dh -e^{-\theta} \int_{-\infty}^{\theta} e^hv(1, h)dh}\\
& = &  \abs{e^{-\theta} \int_{-\infty}^{\theta} e^hv(r(\theta), h)dh -u(1, \theta')+u(1, \theta)-e^{-\theta} \int_{-\infty}^{\theta} e^hv(1, h)dh}\\
& = & \abs{e^{-\theta} \int_{-\infty}^{\theta} e^h [v(r(\theta), h)-v(1, h)]dh -u(1, \theta')+u(1, \theta)}
\end{eqnarray*}
Because of the continuity of $u(1, \theta)$ that we have just shown, there exists a $\delta_3' $ such that 
\begin{equation*}
u(1, \theta)-u(1, \theta') < \frac{\epsilon}{2}, \mbox{ if } \abs{\theta-\theta'} < \delta_3'.
\end{equation*}
According to our assumption (equation (12)), when $\theta$ goes to negative infinity, we have
\begin{equation*}
\abs{v(r(\theta), \theta)-v(1, \theta)} = 0.
\end{equation*}
Therefore, there exists a $\delta_3 '' $ such that
\begin{equation*}
\abs{e^{-\theta} \int_{-\infty}^{\theta} e^h [v(r(\theta), h)-v(1, h)]dh} <  \frac{\epsilon}{2},
\end{equation*}
if $\theta \rightarrow -\infty$ which is the same with $ \abs{r(\theta) -1}<\delta_3''$.

Let $\delta_{3a}=\min{\delta_3', \delta_3''}$. If $\abs{\theta-\theta'} < \delta_{3a}$ and $\abs{r(\theta) -r_0}<\delta_{3a}$, then we have
\begin{equation*}
\abs{ u(r(\theta), \theta)-u(1, \theta')} =\abs{e^{-\theta} \int_{-\infty}^{\theta} e^h [v(r(\theta), h)-v(1, h)]dh -u(1, \theta')+u(1, \theta)} < \epsilon
\end{equation*}
This proves the continuity between $u(r(\theta), \theta)$ and $u(1, \theta)$ when $\theta \rightarrow -\infty$.

Secondly, we show the continuity between $u(r(\theta), \theta)$ and $u(2, \theta)$ when $\theta \rightarrow +\infty$.
\begin{eqnarray*}
& &\abs{ u(r(\theta), \theta)-u(2, \theta')}  \\
& = & \abs{e^{-\theta} \int_{-\infty}^{\theta} e^hv(r(\theta), h)dh-e^{-\theta'}\int_{-\infty}^{\theta'} e^hv(2, h)dh} \\
& = &  \abs{e^{-\theta} \int_{-\infty}^{\theta} e^hv(r(\theta), h)dh-e^{-\theta'}\int_{-\infty}^{\theta'} e^hv(2, h)dh\\
& &+ e^{-\theta} \int_{-\infty}^{\theta} e^hv(2, h)dh -e^{-\theta} \int_{-\infty}^{\theta} e^hv(2, h)dh}\\
& = &  \abs{e^{-\theta} \int_{-\infty}^{\theta} e^hv(r(\theta), h)dh -u(2, \theta')+u(2, \theta)-e^{-\theta} \int_{-\infty}^{\theta} e^hv(2, h)dh}\\
& = & \abs{e^{-\theta} \int_{-\infty}^{\theta} e^h [v(r(\theta), h)-v(2, h)]dh -u(2, \theta')+u(2, \theta)}
\end{eqnarray*}
Because of the continuity of $u(2, \theta)$ that we have just shown, there exists a $\delta_3' $ such that 
\begin{equation*}
u(2, \theta)-u(2, \theta') < \frac{\epsilon}{2}, \mbox{ if } \abs{\theta-\theta'} < \delta_3'.
\end{equation*}
According to our assumption (equation (13)), when $\theta$ goes to positive infinity, we have
\begin{equation*}
\abs{v(r(\theta), \theta)-v(2, \theta)} = 0.
\end{equation*}
Therefore, there exists a $\delta_3 '' $ such that
\begin{equation*}
\abs{e^{-\theta} \int_{-\infty}^{\theta} e^h [v(r(\theta), h)-v(2, h)]dh} <  \frac{\epsilon}{2},
\end{equation*}
if $\theta \rightarrow \infty$ which is the same with $ \abs{r(\theta) -2}<\delta_3''$.

Let $\delta_{3b}=\min{\delta_3', \delta_3''}$. If $\abs{\theta-\theta'} < \delta_{3b}$ and $\abs{r(\theta) -r_0}<\delta_{3b}$, then we have
\begin{equation*}
\abs{ u(r(\theta), \theta)-u(2, \theta')} =\abs{e^{-\theta} \int_{-\infty}^{\theta} e^h [v(r(\theta), h)-v(2, h)]dh -u(2, \theta')+u(2, \theta)} < \epsilon
\end{equation*}
This proves the continuity between $u(r(\theta), \theta)$ and $u(2, \theta)$ when $\theta \rightarrow +\infty$.

\subsection{Annulus with spirals}

Now we change the problem setting again. We still have two concentric circles with radius 1 and 2 respectively, and infinitely many spirals between them. The spirals go though every point between the two circles and each of them goes to approximate both of the circles (Shown as Figure 4). 
\begin{figure}[h]
	\centering
		\includegraphics[width=0.40\textwidth]{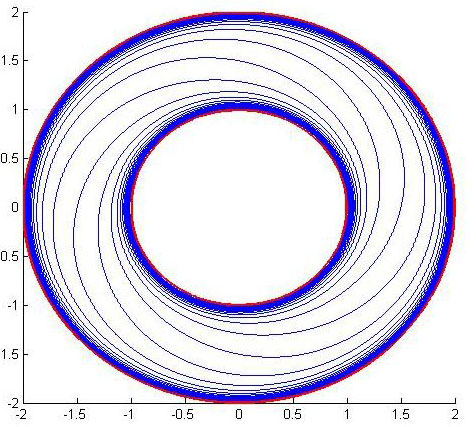}
	\label{fig:figure4}
\end{figure}

The spirals  will neither touch each other nor the two circles so that if we pick an arbitrary small part of the foliation, we can stretch the curves to be straight parallel lines. A continuous and bounded function $V(x, y)$ is defined on the two circles and the spiral, i.e. it is defined on the closed area between the two circles. We want to find another continuous and bounded  function $U(x, y)$ which is differentiable along the spirals and the circles such that $V(x, y)=\nabla U(x, y)$ where $\nabla$ is the gradient of function $U(x, y)$. $\nabla$ has the form of
\begin{equation}
\nabla = F(x, y) \frac{\partial}{\partial x}+ G(x, y) \frac{\partial}{\partial y},
\end{equation}
where $F(x, y)$ and $G(x, y)$ are continuous functions of $x$ and $y$.

Since the definition of the functions are on circles and spirals, it would be more convenient if we transform the coordinates into polar systems. 

Let $r$, $\theta$ be two variables in the corresponding polar system. Then we have
\begin{align}
 x & =r \cos \theta \\
 y & =r \sin \theta
\end{align}
Notice that the radius $r$ is a function of the angle $\theta$, we have the expression of $r$ on spiral s:
\begin{equation}
r(\theta, s)=\frac{3}{2}+\arctan (\theta+s), 
\end{equation}
where $ -\pi \leq s \leq  \pi$. 

We also need the expression for $r$ and $\theta$ in terms of $x$ and $y$ so that we know which spiral is point $(x, y)$ on.
From equation (22), 
\begin{equation}
\theta +s =\tan [\pi (\sqrt{x^2+y^2}-\frac{3}{2})]
\end{equation}
From (20) (21), 
\begin{equation}
\theta=\arctan \frac{y}{x}
\end{equation}
Therefore, we can solve for the expressions of $\theta$ and $s$ in terms of $x$ and $y$.
\begin{align}
 \theta &=\arctan \frac{y}{x} \\
 s &=\tan [\pi (\sqrt{x^2+y^2}-\frac{3}{2})]-\arctan \frac{y}{x}
\end{align}
Notice that the expressions for both $\theta$ and $s$ are uniformly continuous. That indicates that if $V(x, y)$ is continuous, then $v(r, \theta)$ is continuous and if $U(x, y)$ is continuous, then $U(r, \theta)$ is continuous and vice versa.

So we know the transformation from Cartesian to Polar coordinates is:
\begin{align*}
& V(x, y)=V(r \cos \theta, r \sin \theta) = v(r, \theta) \\
& U(x, y)=U(r \cos \theta, r \sin \theta) = u(r, \theta) 
\end{align*} 
The transformation from Polar to Cartesian coordinates is:
\begin{align*}
& v(r, \theta)=v(\sqrt{x^2+y^2}, \arctan \frac{y}{x}) = V(x, y) \\
& u(r, \theta)=u(\sqrt{x^2+y^2}, \arctan \frac{y}{x}) = U(x, y)
\end{align*} 
Then, we do the transformation of gradient. 

The gradient of function $u(r(\theta, s), \theta)$ is
\begin{equation*}
\tilde{\nabla} = \frac{dr}{d\theta} \frac{\partial }{\partial r}+\frac{\partial }{\partial \theta}
\end{equation*}
Fix $s$, from (21) (22) we get
\begin{align*}
 x_s(\theta) & =r \cos \theta =r(\theta, s) \cos \theta\\
 y_s(\theta) & =r \sin \theta =r(\theta, s) \sin \theta
\end{align*}
\begin{eqnarray*}
 \nabla & = & \left[ \begin{array}{cc}
\frac{\partial x}{\partial r}  & \frac{\partial x}{\partial \theta} \\
\frac{\partial y}{\partial r} & \frac{\partial y}{\partial \theta}  \end{array} \right] \tilde{\nabla}\\
& = &  \left[ \begin{array}{cc}
\cos \theta  & -r \sin \theta \\
\sin \theta & r \cos \theta  \end{array} \right] 
 \left[ \begin{array}{c}
\frac{dr}{d\theta} \\
1  \end{array} \right] \\
& = &  \left[ \begin{array}{c}
\cos \theta \frac{dr}{d\theta} -r \sin \theta\\
\sin \theta \frac{dr}{d\theta} +r \cos \theta \end{array} \right] \\
& = & \left[ \begin{array}{c}
F(x, y)\\
G(x, y) \end{array} \right] 
\end{eqnarray*}
Hence we know
\begin{align*}
F(x, y) &= \cos \theta \frac{dr}{d\theta} -r \sin \theta\\
G(x, y)&= \sin \theta \frac{dr}{d\theta} +r \cos \theta\\
V(x, y) &= (F(x, y)\frac{\partial}{\partial x}+G(x, y)\frac{\partial}{\partial y})U(x, y)
\end{align*}

Then, we have to make sure that the transformed variables $x$ and $y$ are eligable under our conditions. 
By assumption, $-\pi < s < \pi$ and $-\infty < \theta < +\infty$. 
It is easy to see that $\arctan \frac{y}{x}$ goes from $-\infty$ to $+\infty$. However, the value of $\arctan \frac{y}{x}$ gives infinite many choices for $\theta$. So we look at equation (26). Given $-\pi < s < \pi$, $\theta$ can only be in an interval with a length smaller than $2\pi$. According to the properties of $\arctan $, we have
\begin{align*}
& \theta = \arctan \frac{y}{x}, \mbox{ if } 2n\pi-\frac{\pi}{2} \leq \theta \leq 2n\pi+\frac{\pi}{2}\\
& \theta = \arctan \frac{y}{x}+\pi, \mbox{ if } 2n\pi+\frac{\pi}{2} \leq \theta \leq 2n\pi+\frac{3\pi}{2}
\end{align*}
where $n$ is a integer.
Thus we only have a unique value for $\theta$ and this transformation is eligable.

Now we start to investigate the properties of all the curves.

On the circles, we have
\begin{equation*}
v(1, \theta)=v(1, \theta+2\pi),
\end{equation*}
\begin{equation*}
v(2, \theta)=v(2, \theta+2\pi),
\end{equation*}
\begin{equation*}
u(1, \theta) = u(1, \theta +2\pi), 
\end{equation*}
\begin{equation*}
u(2, \theta) = u(2, \theta +2\pi).
\end{equation*}
On spiral s, function $v(r(\theta, s), \theta)$ has following property on the spiral:
\begin{equation}
\lim_{\theta \rightarrow -\infty} \abs{ v(r(\theta, s), \theta)-v(1, \theta)} = 0,
\end{equation}
\begin{equation}
\lim_{\theta \rightarrow +\infty} \abs{ v(r(\theta, s), \theta)-v(2, \theta)} = 0.
\end{equation}
We want to show that function $u(r(\theta, s), \theta)$ has the same property:
\begin{equation*}
\lim_{\theta \rightarrow -\infty} \abs{ u(r(\theta, s), \theta)-u(1, \theta)} = 0,
\end{equation*}
\begin{equation*}
\lim_{\theta \rightarrow +\infty} \abs{ u(r(\theta, s), \theta)-u(2, \theta)} = 0.
\end{equation*}

If we set $v(r(\theta, s), \theta)=1$, then $du(r(\theta, s), \theta)/d\theta=1 $ and hence $u(r(\theta, s), \theta) = \theta +C$ which is not a bounded function. This indicates that the differential equation $du(r(\theta, s), \theta)/d\theta=v(r(\theta, s), \theta)$ is not solvable. Therefore, $V(x, y)=\nabla U(x, y)$ is not solvable either. To solve the problem, we add $u(r(\theta), \theta)$ to the left side of the equation and get a new equation $u(r(\theta, s), \theta) + du(r(\theta, s), \theta)/d\theta=v(r(\theta, s), \theta)$. We found that we can solve this new equation by multiplying $e^\theta$ to both sides of the equation and get,
\begin{equation*}
e^\theta u(r(\theta, s), \theta) + e^\theta \frac{d}{d\theta}u(r(\theta, s), \theta)=e^\theta v(r(\theta, s), \theta),
\end{equation*} 
which is the same as 
\begin{equation*}
 \frac{d}{d\theta}e^\theta u(r(\theta, s), \theta))=e^\theta v(r(\theta, s), \theta).
 \end{equation*}
 So, our original equation based on the Cartesian coordinates becomes $U(x, y)+\nabla U(x, y)=V(x, y)$.
 
 Now we are going to find the general form for $U(x, y)$ through finding solution for $u(r(\theta, s), \theta)$.
On the circles, we have
\begin{equation*}
\frac{d}{d\theta}e^{\theta} u(r_0, \theta) \mbox{,  }r_0=1,2.
\end{equation*}
\begin{equation*}
u(r_0, \theta)=e^{-\theta} \int_{-\infty}^{\theta} e^hv(r_0, h)dh \mbox{, } r_0=1,2.
\end{equation*}
Fix a point on a circle with coordinates $(r_0, \theta_0)$, according to the formulae we derived above, we have
\begin{equation}
u(r_0, \theta_0)=e^{-\theta_0} \int_{-\infty}^{\theta_0} e^hv(r_0, h)dh \mbox{, r=1,2.}
\end{equation}
Then, we need to make sure that the function $u(r_0, \theta) $ is periodic with period of $2\pi$, i.e. we need to verify that the equation $u(r_0, \theta_0+2\pi) = u(r_0, \theta_0)$ holds.\\
From equation (23),
\begin{equation*}
u(r_0, \theta_0+2\pi)=e^{-(\theta_0+2\pi)} \int_{-\infty}^{\theta_0+2\pi} e^hv(r_0, h)dh
\end{equation*}
Let $h' = h - 2\pi$,
\begin{eqnarray*}
u(r_0, \theta_0+2\pi) & = & e^{-(\theta_0+2\pi)} \int_{-\infty}^{\theta_0} e^{h'+2\pi}v(r_0, h'+2\pi)dh\\
                                      & = & e^{-\theta_0} \int_{-\infty}^{\theta_0} e^{h'}v(r_0, h')dh\\
                                      & = & u(r_0, \theta_0)
\end{eqnarray*}
This proves that $u(r_0, \theta)$ is a periodic function with period of $2\pi$.\\
Now we will get the solution formulae on the spirals.
Let
\begin{equation}
U'(\theta, s) = u(r(\theta, s), \theta)
\end{equation}
\begin{equation}
V'(\theta, s) = v(r(\theta, s), \theta)
\end{equation}
On the spiral, we have,
\begin{equation*}
\frac{d}{d\theta} e^{\theta} U'(\theta, s)=e^{\theta} V'(\theta, s),
\end{equation*}
\begin{equation*}
U'(\theta, s)=e^{-\theta} \int_{-\infty}^{\theta} e^hV'(h, s)dh.
\end{equation*}
Substitute (30), (31) back to the solution, we get
\begin{equation*}
u(r(\theta, s), \theta) = e^{-\theta} \int_{-\infty}^{\theta} e^hv(r(h, s), h)dh.
\end{equation*}
Fix a point on the spiral with coordinates $(r(\theta_0), \theta_0)$, according to the formulae we derived above, we have
\begin{equation}
u(r(\theta_0, s), \theta_0) = e^{-\theta_0} \int_{-\infty}^{\theta_0} e^hv(R(h, s), h)dh.
\end{equation}

Next step is to show the continuity property of the solution. We are going to prove the following theorem:

\begin{thm}
The function $U(x, y)$ is continuous if the function $V(x, y)$ is continuous.
\end{thm}

Now we show the continuity of the solution we got. We are going to examine the continuity of $e^\theta u(r_0, \theta)(r_0=1,2)$, $e^\theta u(r(\theta, s), \theta)$, between the spirals and between the spirals and the circles. We know $v(r, \theta)$ is continuous since $V(x, y)$ is continuous.

Before we show the continuity of all the things listed above, we first examine the continuity of the function $r(\theta, s)$ in case we need it to know the properties of the foliation.

Let $\epsilon$ be an arbitrary small number. We will show that there exist a real number $\delta_r$ such that $\abs{r(\theta, s)-r(\theta', s')} < \epsilon$ holds for every pair of $(\theta, \theta')$ and $(s, s')$ if $\abs{\theta-\theta'} < \delta_r$ and $\abs{s-s'<\delta_r}$. We have
\begin{equation*}
\abs{r(\theta, s)-r(\theta', s')}=\abs{\arctan (\theta+s)-\arctan (\theta' +s')}
\end{equation*}
Because of the continuity of the $\arctan$ function, there exists a $\delta_r$ such that
\begin{equation*}
\abs{r(\theta, s)-r(\theta', s')} < \epsilon,
\end{equation*}
if $\abs{\theta-\theta'} < \delta_r$ and $\abs{s-s'<\delta_r}$.

\bf Step 1. \normalfont

To show the continuity of the whole solution, we start with showing that the solution on the circles is continuous.
Let $\epsilon$ be an arbitrary small number. We will show that there exist a real number $\delta_1$ such that $\abs{e^\theta u(r_0, \theta)-e^{\theta'} u(r_0, \theta')} < \epsilon$ holds for every pair of $(\theta, \theta')$ if $\abs{\theta-\theta'} < \delta_1$.
Without loss of generality, we assume $\theta>\theta'$.
\begin{eqnarray*}
& &\abs{e^\theta u(r_0, \theta)-e^{\theta'} u(r_0, \theta')} ,     (r_0=1,2)\\
& = & \abs{ \int_{-\infty}^{\theta} e^hv(r_0, h)dh-\int_{-\infty}^{\theta'} e^hv(r_0, h)dh}\\
& = & \abs{\int_{\theta'}^{\theta} e^hv(r_0, h)dh}
\end{eqnarray*}
Because of the continuity of $v(r_0, h)$, there exists a $\delta_1$ such that 
\begin{equation*}
\abs{e^\theta u(r_0, \theta)-e^{\theta'} u(r_0, \theta')} < \epsilon, (r_0=1,2)
\end{equation*}

\bf Step 2.\normalfont

Showing the solution is continuous on each single spiral.
We randomly choose a spiral s with the radius function $ r(\theta, s)$. Since $\theta$ is the only variable here, we can write the radius number of spiral s to be $r_s(\theta)$ which is obviously a continuous function. 

Since both $r_s(h)$ and $v(r, h)$ are continuous, $v(r_s(h), h)$ is again a continuous function.

Let $\epsilon$ be an arbitrary small number. We will show that there exist a real number $\delta_2$ such that $\abs{e^\theta u(r_s(\theta), \theta)-e^{\theta'} u(r_s(\theta'), \theta')} < \epsilon$ holds for every pair of $(\theta, \theta')$ if $\abs{\theta-\theta'} < \delta_2$.
Without loss of generality, we assume $\theta>\theta'$.
\begin{eqnarray*}
& &\abs{e^\theta u(r_s(\theta), \theta)-e^{\theta'} u(r_s(\theta'), \theta')} \\
& = & \abs{ \int_{-\infty}^{\theta} e^hv(r_s(h), h)dh-\int_{-\infty}^{\theta'} e^hv(r_s(h), h)dh}\\
& = & \abs{\int_{\theta'}^{\theta} e^hv(r_s(h), h)dh}
\end{eqnarray*}
Because of the continuity of $v(r_s(h), h)$, there exists a $\delta_2$ such that 
\begin{equation*}
\abs{e^\theta u(r_s(\theta), \theta)-e^{\theta'} u(r_s(\theta'), \theta')} < \epsilon
\end{equation*}

\bf Step 3.\normalfont

Showing the solution is continuous between different spirals. 

Let $\epsilon$ be an arbitrary small number. We will show that there exist a real number $\delta_3$ such that $\abs{e^\theta u(r(\theta, s), \theta)-e^{\theta'} u(r(\theta', s'), \theta')} < \epsilon$ holds for every pair of $(\theta, \theta')$ and $(s, s')$ if $\abs{\theta-\theta'} < \delta_3$ and $\abs{s-s'} < \delta_3$.
Without loss of generality, we assume $\theta>\theta'$. Let $x  =r \cos \theta$, $y  =r \sin \theta$, $x'  =r' \cos \theta'$ and $y'  =r' \sin \theta'$.
\begin{eqnarray*}
& &U(x, y)-U(x', y')\\
& = & \abs{e^\theta u(r(\theta, s), \theta)-e^{\theta'} u(r(\theta', s'), \theta')} \\
& = & \abs{ \int_{-\infty}^{\theta} e^hv(r(h, s), h)dh-\int_{-\infty}^{\theta'} e^hv(r(h, s'), h)dh}\\
& = & \abs{ \int_{-\infty}^{\theta} e^h [v(r(h, s), h)-v(r(h, s'), h)] dh-\int_{\theta}^{\theta'} e^hv(r(h, s'), h)dh}  \\
& =&  \abs{ \int_{-\infty}^{-L} e^h [v(r(h, s), h)-v(r(h, s'), h)] dh+ \int_{-L}^{\theta} e^h [v(r(h, s), h)-v(r(h, s'), h)] dh\\
& & -\int_{\theta}^{\theta'} e^hv(r(h, s'), h)dh}
\end{eqnarray*}
where L is a large positive number to be chosen.

Because of the continuity of $v(r(h, s'), h)$, there exists a $\delta_{3a}$ such that
\begin{equation*}
\int_{\theta}^{\theta'} e^hv(r(h, s'), h)dh  < \frac{\epsilon}{2},
\end{equation*}
if $\abs{\theta-\theta'} < \delta_{3a}$.\\
If $|x-x'|< \delta_x$ and $|y-y'|<\delta_y$, because of the continuity of $s(x, y)$, $|s-s'|$ can be arbitrary small. Then, because of the continuity of $v(r(h, s'), h)$, there exists a small $\delta_{3b}$ such that
\begin{equation*}
v(r(h, s), h)-v(r(h, s'), h) < \frac{\epsilon}{4},
\end{equation*}
if $|s-s'|<\delta_{3b}$ which is the same with $|x-x'|< \delta_x$ and $|y-y'|<\delta_y$.

After that, we evaluate L.

Let 
\begin{equation*}
\abs{ v(r(h, s), h)-v(r(h, s'), h)} \leq M,
\end{equation*} 
where $M$ is a positive real number.

Then we have
\begin{align*}
&\abs{ \int_{-\infty}^{-L} e^h M dh} = \frac{\epsilon}{4}\\
&e^{-L} M \leq  \frac{\epsilon}{4}\\
&L \geq \ln \frac{4 M}{\epsilon}
\end{align*}
Therefore, we have
\begin{equation*}
\abs{ \int_{-\infty}^{-\ln \frac{4 M}{\epsilon}} e^h [v(r(h, s), h)-v(r(h, s'), h)] dh} < \frac{\epsilon}{4},
\end{equation*}
if $|s-s'|<\delta_{3b}$.

Since $v(r(h, s), h)$ is a continuous and bounded function, we have 
\begin{equation*}
v(r(h, s), h)-v(r(h, s'), h) < \frac{\epsilon}{4(\theta - (-\ln \frac{4M}{\epsilon}))e^\theta}
\end{equation*}
if $\abs{s-s'} < \delta_{3b}$ which is the same with $|x-x'|< \delta_x$ and $|y-y'|<\delta_y$

Therefore,
\begin{eqnarray*}
 & & \int_{-\ln \frac{4M}{\epsilon}}^{\theta} e^h [v(r(h, s), h)-v(r(h, s'), h)] dh\\
 &<& \int_{-\ln \frac{4M}{\epsilon}}^{\theta} e^\theta \frac{\epsilon}{4(\theta - (-\ln \frac{4M}{\epsilon}))e^\theta))e^\theta}
 dh \\
 &=& \frac{\epsilon}{4},
 \end{eqnarray*}
 if $\abs{s-s'} < \delta_{3b}$.

Then, let $\delta_3=\min \{\delta_{3a}, \delta_{3b}\} $, we have
\begin{equation*}
\abs{e^\theta u(r(\theta, s), \theta)-e^{\theta'} u(r(\theta', s'), \theta')} < \epsilon,
\end{equation*}
which gives,
\begin{equation*}
|U(x, y)-U(x', y')|<\epsilon,
\end{equation*}
if  $\abs{\theta-\theta'} < \delta_3$ and $\abs{s-s'} < \delta_3$ which is the same with $|x-x'|< \delta_x$ and $|y-y'|<\delta_y$

\bf Step 4.\normalfont

Showing the solution is continuous between each spiral and the inner circle.

We know that $v(r(h, s), h)$ is continuous function for all s.

Let $\epsilon$ be an arbitrary small number. We will show that there exist a real number $\delta_4$ such that $\abs{[u(r(\theta, s), \theta)- u(1, \theta')]} < \epsilon$ holds for every pair of $(\theta, \theta')$ and $(r(\theta, s), 1)$ if $\abs{\theta+ 2n \pi-\theta'} < \delta_4$ and $ \abs{r(\theta, s) -1}<\delta_4$.
Without loss of generality, we assume $\theta > \theta'$. We will start with the cartesian coordinates.
\begin{eqnarray*}
&  &\abs{ u(r(\theta, s), \theta)-u(1, \theta')}  \\
& = & \abs{e^{-\theta} \int_{-\infty}^{\theta} e^hv(r(\theta, s), h)dh-e^{-\theta'}\int_{-\infty}^{\theta'} e^hv(1, h)dh} \\
& = &  \abs{e^{-\theta} \int_{-\infty}^{\theta} e^hv(r(\theta, s), h)dh-e^{-\theta'}\int_{-\infty}^{\theta'} e^hv(1, h)dh\\
& &+ e^{-\theta} \int_{-\infty}^{\theta} e^hv(1, h)dh -e^{-\theta} \int_{-\infty}^{\theta} e^hv(1, h)dh}\\
& = &  \abs{e^{-\theta} \int_{-\infty}^{\theta} e^hv(r(\theta, s), h)dh -u(1, \theta')+u(1, \theta)-e^{-\theta} \int_{-\infty}^{\theta} e^hv(1, h)dh}\\
& = & \abs{e^{-\theta} \int_{-\infty}^{\theta} e^h [v(r(\theta, s), h)-v(1, h)]dh -u(1, \theta')+u(1, \theta)}
\end{eqnarray*}
Because of the continuity of $u(1, \theta)$ that we have just shown, there exists a $\delta_{4a} $ such that 
\begin{equation*}
u(1, \theta)-u(1, \theta') < \frac{\epsilon}{2}, \mbox{ if } \abs{\theta+ 2n \pi-\theta'} < \delta_{4a}.
\end{equation*}
According to our assumption (equation (21)), when $\theta \rightarrow -\infty$, we have
\begin{equation*}
\abs{v(r(\theta, s), \theta)-v(1, \theta)} = 0, \mbox{ for }s\mbox{ uniformly.}
\end{equation*}
Therefore, there exists a $\delta_{4b}$ such that
\begin{equation*}
\abs{e^{-\theta} \int_{-\infty}^{\theta} e^h [v(r(\theta, s), h)-v(1, h)]dh} <  \frac{\epsilon}{2},
\end{equation*}
if $\theta \rightarrow -\infty$ which is the same with $ \abs{r(\theta, s) -1}<\delta_{4b}$.

Let $\delta_{4}=\min{\delta_{4a}, \delta_{4b}}$. If $\abs{\theta-\theta'} < \delta_{4}$ and $\abs{r(\theta, s) -1}<\delta_{4}$, then we have
\begin{equation*}
\abs{ u(r(\theta, s), \theta)-u(1, \theta')} =\abs{e^{-\theta} \int_{-\infty}^{\theta} e^h [v(r_s(\theta), h)-v(1, h)]dh -u(1, \theta')+u(1, \theta)} < \epsilon
\end{equation*}
This proves the continuity between $u(r(\theta, s), \theta)$ and $u(1, \theta)$ for all s when $\theta \rightarrow -\infty$.

\bf Step 5.\normalfont

Showing the solution is continuous between each spiral and the outer circle.

We know that $v(r(s, h), h)$ is continuous function for all s.

Let $\epsilon$ be an arbitrary small number. We will show that there exist a real number $\delta_5$ such that $\abs{[u(r(\theta, s), \theta)- u(2, \theta')]} < \epsilon$ holds for every pair of $(\theta, \theta')$ and $(r(\theta, s), 2)$ when $\theta \rightarrow +\infty$ if $\abs{\theta+2n \pi-\theta'} < \delta_5$ and $ \abs{r(\theta, s) -1}<\delta_5$.
Without loss of generality, we assume $\theta > \theta'$.
\begin{eqnarray*}
& &\abs{ u(r(\theta, s), \theta)-u(2, \theta')}  \\
& = & \abs{e^{-\theta} \int_{-\infty}^{\theta} e^hv(r(\theta, s), h)dh-e^{-\theta'}\int_{-\infty}^{\theta'} e^hv(2, h)dh} \\
& = &  \abs{e^{-\theta} \int_{-\infty}^{\theta} e^hv(r(\theta, s), h)dh-e^{-\theta'}\int_{-\infty}^{\theta'} e^hv(2, h)dh\\
& &+ e^{-\theta} \int_{-\infty}^{\theta} e^hv(2, h)dh -e^{-\theta} \int_{-\infty}^{\theta} e^hv(2, h)dh}\\
& = &  \abs{e^{-\theta} \int_{-\infty}^{\theta} e^hv(r(\theta, s), h)dh -u(2, \theta')+u(2, \theta)-e^{-\theta} \int_{-\infty}^{\theta} e^hv(2, h)dh}\\
& = & \abs{e^{-\theta} \int_{-\infty}^{\theta} e^h [v(r(\theta, s), h)-v(2, h)]dh -u(2, \theta')+u(2, \theta)}
\end{eqnarray*}
Because of the continuity of $u(2, \theta)$ that we have just shown, there exists a $\delta_{5a} $ such that 
\begin{equation*}
u(2, \theta)-u(2, \theta') < \frac{\epsilon}{2}, \mbox{ if } \abs{\theta+2n \pi-\theta'} < \delta_{5a}.
\end{equation*}
According to our assumption (equation (22)), when $\theta \rightarrow +\infty$, we have
\begin{equation*}
\abs{v(r(\theta, s), \theta)-v(2, \theta)} = 0, \mbox{ for }s\mbox{ uniformly.}
\end{equation*}
Therefore, there exists a $\delta_{5b}$ such that
\begin{equation*}
\abs{e^{-\theta} \int_{-\infty}^{\theta} e^h [v(r(\theta, s), h)-v(2, h)]dh} <  \frac{\epsilon}{2},
\end{equation*}
if $\theta \rightarrow +\infty$ which is the same with $ \abs{(\theta, s)-2}<\delta_{4b}$.

Let $\delta_{5}=\min{\delta_{5a}, \delta_{5b}}$. If $\abs{\theta-\theta'} < \delta_{5}$ and $\abs{r_s(\theta) -1}<\delta_{5}$, then we have
\begin{equation*}
\abs{ u(r(\theta, s), \theta)-u(2, \theta')} =\abs{e^{-\theta} \int_{-\infty}^{\theta} e^h [v(r(\theta, s), h)-v(2, h)]dh -u(2, \theta')+u(2, \theta)} < \epsilon
\end{equation*}
This proves the continuity between $u(r(\theta, s), \theta)$ and $u(2, \theta)$ for all s when $\theta \rightarrow +\infty$.

In conlusion, let $\delta = \min \{\delta_1, \delta_2, \delta_3, \delta_4, \delta_5\}$. Then we have 
\begin{equation*}
\abs{ u(r, \theta)-u(r', \theta')} < \epsilon
\end{equation*}
if $|r-r'|<\delta$ and $|\theta-\theta'|<\delta$.

Therefore, $u(r, \theta)$ is continuous and hence $U(x, y)$ is also continuous.

\subsection{General foliations}

Consider a vector field $X$ in $\mathbb{R}^n$, we have a compact set which is a union of integral curves of vector field $X$. So, in particular we assume the the vector field $X$ does not vanish at any point in $K$. Let $V$ be a continuous function on a compact set $K$. Based on previous experience, we want to find a continuous function $U$ which is differentiable along the integral curves so that on each curve,  
\begin{equation}
U+XU=V.
\end{equation}
Imagine there is a particle moving along the curves, $X$ can be viewed as the velocity of the particle and $\norm{X}$ is the speed. We could set $\norm{X}=1$ by changing the coordinates so that the particles move at speed of 1 in all directions on all curves.

Here is another way to state this problem. Consider a compact set $K$ in $\mathbb{R}^n$. There is only one curve through every point in $K$. The curves vary continuously, have similar slopes and do not intersect with each other. Let $V$ be a continuous function on $K$, we will find a continuous function $U$ which is differentiable along the curves so that on each curve, 
\begin{equation}
U+A\frac{dU}{dt}=V,
\end{equation}
where $t$ is the distance that a point on a curve go through and $A$ is a coefficient function which is continuous and positive. Notice that since function $A$ is defined on a compact set $K$, $A$ is bounded automatically.
 
To relate the two case, let $A=\norm{X}=1$, then $X=A\frac{X}{\norm{X}}$. From (34), we have
\begin{equation*}
U+A\frac{X}{\norm{X}}U=U+\frac{d}{dt}U=V,
\end{equation*}
which is the same with (35). Thus we know the two conditions are identical and we can solve the general foliations through solving $U+\frac{dU}{dt}=V$.

We first focus on the real line case. On a real line, our differential equation becomes $U(x)+\frac{dU(x)}{dx}=V(x)$.

Based on previous experience, we multiply $e^x$ to the equation and get:
\begin{align*}
& e^xU(x)+e^x\frac{dU(x)}{dx}=e^xV(x)\\
& \frac{d(e^xU(x)}{dx}=e^xV(x)\\
& e^xU(x)=\int e^xV(x)
\end{align*} 
Let $V(x)=1$, we have $U(x)=1$ which is obviously a bounded function. Hence the equation $U(x)+\frac{dU(x)}{dx}=V(x)$ is solvable.

So we have the general formulae for U(x):
\begin{equation}
U(x)=e^{-x} \int_{-\infty}^{x} e^tV(t)dt.
\end{equation}
If $V(x)$ is periodic functions with period of 1, we want to know whether $U(x)$ is also periodic or not.
Let $t'=t-1$.
\begin{eqnarray*}
U(x+1)&=&e^{-(x+1)} \int_{-\infty}^{x+1} e^tV(t)dt\\
&=&e^{-(x+1)} \int_{-\infty}^{x} e^{t'+1}V(t'+1)dt'\\
&=&e^{-x} \int_{-\infty}^{x} e^{t'}V(t')dt'\\
&=&U(x).
\end{eqnarray*}
So $U(x)$ is a periodic function with period of 1.
 
We also consider another solvable equation $\frac{d(e^xA(x)U(x))}{dx}=e^xV(x)$.
\begin{align*}
& \frac{d(e^xA(x)U(x))}{dx}=e^xV(x)\\
& \frac{dA(x)}{dx}U(x)e^x+\frac{dU(x)}{dx}A(x)e^x+A(x)U(x)e^x=e^xV(x)\\
& \frac{dA(x)}{dx}U(x)+\frac{dU(x)}{dx}A(x)+A(x)U(x)=V(x)\\
& (\frac{dA(x)}{dx}+A(x))U(x)+A(x)\frac{dU(x)}{dx}=V(x)
\end{align*}
Set V(x)=1, we have
\begin{align*}
&A(x)U(x)e^x=\int 1e^xdx\\
&U(x)=\frac{1}{A(x)}
\end{align*}
which is obviously a bounded function. 

Then we know our modified equation $\frac{d(e^xA(x)U(x))}{dx}=e^xV(x)$ which is the same with $(\frac{dA(x)}{dx}+A(x))U(x)+A(x)\frac{dU(x)}{dx}=V(x)$ is solvable.

So we have the general formulae for U(x):
\begin{equation}
U(x)=\frac{1}{A(x)e^x} \int_{-\infty}^{x} e^tV(t)dt.
\end{equation}
Going back to the original case, our equation is:
\begin{equation}
(\norm{X}+\frac{d\norm{X}}{ds})U+XU=V.
\end{equation}

If $V(x)$ and $A(x)$ are periodic functions with periods of 1, we want to know whether $U(x)$ is also periodic or not.
Let $t'=t-1$.
\begin{eqnarray*}
U(x+1)&=&\frac{1}{A(x+1)e^{x+1}} \int_{-\infty}^{x+1} e^tV(t)dt\\
&=&\frac{1}{A(x)e^{x+1}} \int_{-\infty}^{x} e^{t'+1}V(t'+1)dt'\\
&=&\frac{1}{A(x)e^x} \int_{-\infty}^{x} e^{t'}V(t')dt'\\
&=&U(x).
\end{eqnarray*}
So $U(x)$ is a periodic function with period of 1. 
But we are not going to discuss more about this equation in this paper.

We have the following assumptions:
Let vector field $X \in C^\infty$, and $X \neq 0$.
Let $K$ be a compact union of integral curves.
Then there exists a map $\Phi$
\begin{equation}
\Phi : \mathbb{R}^k \times \mathbb{R} \rightarrow \mathbb{R}^k
\end{equation}
Let $x=(x_1, x_2, ... , x_k) \in \mathbb{R}^k$ be an initial point on the compact set $K$. Then, by assumption, there is only one curve $\gamma_x(t)$ goes thought $x$ parameterized by time $t$. When $t=0$, $\gamma_x(0)=x$ since $x$ is the initial point. Then, we map the initial point to another point after time $t$ by $\Phi$. We have
\begin{equation}
\Phi (x, t) = \gamma_x(t).
\end{equation}
Since $X \in C^\infty$, the map $\Phi \in C^\infty$.
For each $t_0 \in \mathbb{R}$, the map $x \rightarrow \Phi (x, t_0)$ has a $C^\infty$ inverse map $\phi_{t_0}$ in a neighborhood of $K$.
\begin{equation}
\phi_{t_0}(\Phi (x, t_0))=x
\end{equation}

\begin{thm}
On a compact set $K$, if $V$ is continuous along an integral curve, then $U$ is continuous along the curve.
\end{thm}

Proof. Pick any integral curve $\gamma_x(t)$, we have a function $v(t)$ which is defined as
\begin{equation}
v(t)=V(\gamma_x(t))
\end{equation}
Then we have a differential equation $u(t)+\frac{d}{dt}u(t)=v(t)$ where $u(t)$ is a bounded function which is differentiable along the curve. We know that
\begin{equation}
u(t)=U(\gamma_x(t))
\end{equation}
So what we need to show is the continuity of $u(t)$ given that $v(t)$ is continuous.
We have the general formulae of $u(t)$ which is:
\begin{equation*}
u(t)=e^{-t} \int_{-\infty}^{t} e^wv(w)dw
\end{equation*}
We will prove the continuity of $u(t)$ by showing that $e^tu(t)$ is continuous. 

Let $\epsilon$ be an arbitrary small number. We will show that there exist a real number $\delta$ such that $\abs{e^tu(t)-e^{t'}u(t')}<\epsilon$ holds for every pair of $(t, t')$ if $\abs{t-t'}<\delta$. Without loss of generality, we assume $t>t'$.
\begin{eqnarray*}
&&\abs{e^tu(t)-e^{t'}u(t')}\\
&=&\abs{\int_{-\infty}^t e^wv(w)dw -\int_{-\infty}^{t'} e^wv(w)}\\
&=&\abs{\int_{t'}^t e^wv(w)dw}
\end{eqnarray*}
Let $k=\max \{e^tv(t), e^{t'}v(t')\} $, because of the continuity of $v(t)$, there exist a $\delta=|t-t'|$ such that
\begin{equation*}
\abs{\int_{t'}^t e^wv(w)dw}<(t-t')k=\delta k <\epsilon
\end{equation*}
\begin{equation*}
\delta <\frac{\epsilon}{k}
\end{equation*}
If $|t-t'|<\delta$, then $\abs{e^tu(t)-e^{t'}u(t')}<\epsilon$ and hence $u(t)$ is continuous. Therefore, $U(t)$ is continuous. 

\begin{thm}
On a compact set $K$, if $V$ is continuous function, then $U$ is continuous.
\end{thm}

Proof. According to what we have done above, we have
\begin{equation}
U(\gamma_x(t))=e^{-t} \int_{-\infty}^t e^wV(\gamma_x(w)dw
\end{equation}
When $t=0$, we know that $\gamma_x(0)=x$ since $x$ is the initial point. Then we have
\begin{equation}
U(x)=\int_{-\infty}^0 e^wV(\gamma_x(w)dw
\end{equation}
By assumption, we have
\begin{equation}
\Phi (x, t) = \gamma_x(t)
\end{equation}
So our formulae for $U$ when $t=0$ is
\begin{equation}
U(x)=\int_{-\infty}^0 e^wV(\Phi (x, w))dw
\end{equation}
Let $\epsilon$ be an arbitrary small number. We will show that there exist a real number $\delta$ such that $\abs{U(x)-U(x')}<\epsilon$ holds for every pair of $(x, x')$ if $\abs{x-x'}<\delta$. Let $L$ be a large number to be chosen.
\begin{eqnarray*}
&&\abs{U(x)-U(x')}\\
&=&\abs{\int_{-\infty}^0 e^wV(\Phi (x, w))dw-\int_{-\infty}^0 e^wV(\Phi (x', w))dw}\\
&=&\abs{\int_{-\infty}^{-L} [e^wV(\Phi (x, w))- e^wV(\Phi (x', w))]dw+\int_{-L}^0 [e^wV(\Phi (x, w))- e^wV(\Phi (x', w))]dw}
\end{eqnarray*}
Let $M$ be the maximum of function $V$, then
\begin{equation}
V(\Phi(x, w))-V(\Phi(x', w))<2M
\end{equation}
\begin{eqnarray*}
&&\int_{-\infty}^{-L} [e^wV(\Phi (x, w))- e^wV(\Phi (x', w))]dw\\
&<&\int_{-\infty}^{-L} e^w 2Mdw\\
&=&2Me^{-L}<\frac{\epsilon}{2}
\end{eqnarray*}
Solve the inequation above, we get 
\begin{equation}
L>ln \frac{4M}{\epsilon}
\end{equation}
On the other hand, because of the continuity of $V$ and $\Phi$, there exists a $delta$ such that
\begin{equation*}
\int_{-ln \frac{4M}{\epsilon}}^0 [e^wV(\Phi (x, w))- e^wV(\Phi (x', w))]dw<\frac{\epsilon}{2}
\end{equation*}
if $\abs{x-x'}<\delta$. Therefore, $U$ is continuous on $K$.

Then we are going to find that if $V$ is better than continuous, what will happen to $U$. 

\begin{thm}
On a compact set $K$, if  $V \in C^1$, then $U \in C^1$.
\end{thm}

Proof.  If $U \in C^1$, from (47), we will have
\begin{equation*}
\frac{\partial U}{\partial x_i}=\int_{-\infty}^{0} \frac{\partial }{\partial x}(V(\Phi(x, w)))dw
\end{equation*}
We will show that equation above holds.
Let 
\begin{equation}
f(x, w)=V(\Phi(x, w))
\end{equation}
Notice that $f(x) \in C^1$.
Then the equation that we are going to prove becomes
\begin{equation}
\frac{\partial U}{\partial x_i}=\int_{-\infty}^{0} \frac{\partial }{\partial x}f(x, w)dw
\end{equation}
By definition of derivative, we have
\begin{equation*}
\frac{\partial U}{\partial x_i}=\lim_{\Delta x \rightarrow 0} \frac{U(x+\Delta x)-U(x)}{\Delta x}.
\end{equation*}
By the Mean-value theorem,
\begin{eqnarray*}
&&\frac{U(x+\Delta x)-U(x)}{\Delta x}=\int_{-\infty}^0 e^w \frac{f(x+ \Delta x, w)-f(x, w)}{\Delta x}\\
&=&\int_{-\infty}^0 e^w \frac{\partial }{\partial x_i} (f(x'(w), w)dw 
\end{eqnarray*}
where $x \leq x'(w) \leq x+\Delta x$.
\begin{equation*}
\int_{-\infty}^0 e^w \frac{\partial }{\partial x_i} f(x'(w), w)dw =\int_{-\infty}^{-L} e^w \frac{\partial }{\partial x_i} f(x'(w), w)dw +\int_{-L}^0 e^w \frac{\partial }{\partial x_i} f(x'(w), w)dw 
\end{equation*}
\begin{equation*}
\int_{-\infty}^0 e^w \frac{\partial }{\partial x_i}f(x, w)dw =\int_{-\infty}^{-L} e^w \frac{\partial }{\partial x_i}f(x, w)dw +\int_{-L}^0 e^w \frac{\partial }{\partial x_i}f(x, w)dw 
\end{equation*}
where $L$ is a large number to be chosen.

Let $M$ be the maximum value of $\frac{\partial }{\partial x_i}(f(x, w)$ and $M'$ be the maximum value of $\frac{\partial }{\partial x_i}f(x'(w), w)$. Let $\epsilon$ be an arbitrary small number.
\begin{eqnarray*}
&&\int_{-\infty}^{-L} e^w \frac{\partial }{\partial x_i}f(x'(w), w)dw\\
&< &\int_{-\infty}^{-L} e^wM'dw\\
&=& M'e^{-L} < \epsilon
\end{eqnarray*}
Solve the inequation above, we get
\begin{equation}
L > ln \frac{M'}{\epsilon}
\end{equation}
Similarly, we can get 
\begin{equation}
L > ln \frac{M}{\epsilon}
\end{equation}
by setting
\begin{equation*}
\int_{-\infty}^{-L} e^w \frac{\partial }{\partial x_i}f(x, w)dw < \epsilon.
\end{equation*}
Then we make $L > max \{ln \frac{M'}{\epsilon}, ln \frac{M}{\epsilon} \}$ and get the following equations:
\begin{equation*}
\int_{-\infty}^0 e^w \frac{\partial }{\partial x_i}f(x'(w), w)dw = \int_{-L}^0 e^w \frac{\partial }{\partial x_i}f(x'(w), w)dw + \epsilon
\end{equation*}
\begin{equation*}
\int_{-\infty}^0 e^w \frac{\partial }{\partial x_i}f(x, w)dw =\int_{-L}^0 e^w \frac{\partial }{\partial x_i}f(x, w)dw + \epsilon
\end{equation*}
When $\Delta x \rightarrow 0$, $x'(w) \rightarrow x$. 

Hence we know that
\begin{eqnarray*}
\frac{dU}{dx}&=&\lim_{\Delta x \rightarrow 0} \frac{U(x+\Delta x)-U(x)}{\Delta x}\\
&=&\lim_{\Delta x \rightarrow 0} \int_{-\infty}^0 e^w \frac{\partial }{\partial x_i}f(x'(w), w)dw \\
&= &\int_{-\infty}^0 e^w \frac{\partial }{\partial x_i}f(x, w)dw
\end{eqnarray*}
Therefore, $U \in C^1$.

\begin{thm}
On a compact set $K$, if  $V \in C^n$, then $U \in C^n$.
\end{thm}
Proof. We will show this theorem by mathematical induction. 

Initial Step. 
From theorem 4, we know when $n=1$, if $V \in C^1$, then $U \in C^1$.

Induction Step. 
Assume that for $n=k$, if $V \in C^k$, then $U \in C^k$. That is to say that we have
\begin{equation}
\frac{\partial^k U }{\partial x_{i_1} ... \partial x_{i_k}} = \int_{-\infty}^0 e^w \frac{\partial^k}{\partial x_{i_1} ... \partial x_{i_k}} f(x, w)dw
\end{equation}
We will show that for $n=k+1$, if $V \in C^{k+1}$, then $U \in C^{k+1}$. That is to say, we will show that the following equation holds:
\begin{equation}
\frac{\partial^{k+1} U }{\partial x_{i_1} ... \partial x_{i_{k+1}}} = \int_{-\infty}^0 e^w \frac{\partial^{k+1}}{\partial x_{i_1} ... \partial x_{i_{k+1}}} f(x, w)dw
\end{equation}
Let 
\begin{equation*}
\tilde{f}(x, w)= \frac{\partial^k}{\partial x_{i_1} ... \partial x_{i_k}} f(x, w).
\end{equation*}
Since $f(x, w) \in C^{k+1}$, it is obvious that $\tilde{f}(x, w) \in C^1$. 
Let
\begin{equation*}
\tilde{U}= \frac{\partial^k U}{\partial x_{i_1} ... \partial x_{i_k}}.
\end{equation*}
We will show that 
\begin{equation}
\frac{\partial \tilde{U}}{\partial x_{i_{k+1}}}=\int_{-\infty}^{0} \frac{\partial }{\partial x_{i_{k+1}}}\tilde{f}(x, w)))dw
\end{equation}
By definition of derivative, we have
\begin{equation*}
\frac{\partial \tilde{U}}{\partial x_{i_{k+1}}}=\lim_{\Delta x \rightarrow 0} \frac{\tilde{U}(x+\Delta x)-\tilde{U}(x)}{\Delta x}.
\end{equation*}
By the Mean-value theorem,
\begin{eqnarray*}
&&\frac{\tilde{U}(x+\Delta x)-\tilde{U}(x)}{\Delta x}=\int_{-\infty}^0 e^w \frac{\tilde{f}(x+ \Delta x, w)-\tilde{f}(x, w)}{\Delta x}\\
&=&\int_{-\infty}^0 e^w \frac{\partial }{\partial x_{i_{k+1}}} \tilde{f}(x'(w), w)dw 
\end{eqnarray*}
where $x \leq x'(w) \leq x+\Delta x$.
\begin{equation*}
\int_{-\infty}^0 e^w \frac{\partial }{\partial x_{i_{k+1}}} \tilde{f}(x'(w), w)dw =\int_{-\infty}^{-L} e^w \frac{\partial }{\partial x_{i_{k+1}}} \tilde{f}(x'(w), w)dw +\int_{-L}^0 e^w \frac{\partial }{\partial x_{i_{k+1}}} \tilde{f}(x'(w), w)dw 
\end{equation*}
\begin{equation*}
\int_{-\infty}^0 e^w \frac{\partial }{\partial x_{i_{k+1}}}\tilde{f}(x, w)dw =\int_{-\infty}^{-L} e^w \frac{\partial }{\partial x_{i_{k+1}}}\tilde{f}(x, w)dw +\int_{-L}^0 e^w \frac{\partial }{\partial x_{i_{k+1}}}\tilde{f}(x, w)dw 
\end{equation*}
where $L$ is a large number to be chosen.

Let $M$ be the maximum value of $\frac{\partial }{\partial x_{i_{k+1}}}\tilde{f}(x, w)$ and $M'$ be the maximum value of $\frac{\partial }{\partial x_{i_{k+1}}}\tilde{f}(x'(w), w)$. Let $\epsilon$ be an arbitrary small number.
\begin{eqnarray*}
&&\int_{-\infty}^{-L} e^w \frac{\partial }{\partial x_{i_{k+1}}}\tilde{f}(x'(w), w)dw\\
&< &\int_{-\infty}^{-L} e^wM'dw\\
&=& M'e^{-L} < \epsilon
\end{eqnarray*}
Solve the inequation above, we get
\begin{equation}
L > ln \frac{M'}{\epsilon}
\end{equation}
Similarly, we can get 
\begin{equation}
L > ln \frac{M}{\epsilon}
\end{equation}
by setting
\begin{equation*}
\int_{-\infty}^{-L} e^w \frac{\partial }{\partial x_{i_{k+1}}}\tilde{f}(x, w)dw < \epsilon.
\end{equation*}
Then we make $L > max \{ln \frac{M'}{\epsilon}, ln \frac{M}{\epsilon} \}$ and get the following equations:
\begin{equation*}
\int_{-\infty}^0 e^w \frac{\partial }{\partial x_{i_{k+1}}}\tilde{f}(x'(w), w)dw = \int_{-L}^0 e^w \frac{\partial }{\partial x_{i_{k+1}}}\tilde{f}(x'(w), w)dw + \epsilon
\end{equation*}
\begin{equation*}
\int_{-\infty}^0 e^w \frac{\partial }{\partial x_{i_{k+1}}}\tilde{f}(x, w)dw =\int_{-L}^0 e^w \frac{\partial }{\partial x_{i_{k+1}}}\tilde{f}(x, w)dw + \epsilon
\end{equation*}
When $\Delta x \rightarrow 0$, $x'(w) \rightarrow x$. 

Hence we know that
\begin{eqnarray*}
\frac{\partial^{k+1} U }{\partial x_{i_1} ... \partial x_{i_{k+1}}} &=&\frac{\partial \tilde{U}}{\partial x_{i_{k+1}}}\\
&=&\lim_{\Delta x \rightarrow 0} \frac{\tilde{U}(x+\Delta x)-\tilde{U}(x)}{\Delta x}\\
&=&\lim_{\Delta x \rightarrow 0} \int_{-\infty}^0 e^w \frac{\partial }{\partial x_{i_{k+1}}}\tilde{f}(x'(w), w)dw \\
&= &\int_{-\infty}^0 e^w \frac{\partial }{\partial x_{i_{k+1}}}\tilde{f}(x, w)dw\\
&=& \int_{-\infty}^0 e^w \frac{\partial^k}{\partial x_{i_1} ... \partial x_{i_k}} f(x, w)dw
\end{eqnarray*}
Therefore, $\tilde{U} \in C^1$ and  $U \in C^{k+1}$.

\section{Foliations on vector field with singularities}

In the previous section, we only considered nice vector field $X \neq 0$ everywhere. Do we really need that condition? In this section, we are going to check with vector field $X= \sum \phi (x_i)\frac{\partial}{\partial X_i} $ on a compact set $K$ with singularities at the origin and see whether that works well. 

\subsection{The real line case}

Since we still need a bounded vector field, $ \phi(x)= \arctan x$ would be a good choice. However, it would be a huge mass when we do the calculations with $\arctan x$, so we modify function $\phi$ into:
\[ \phi(x) = \left\{ \begin{array}{ll}
         1 & \mbox{if $x \geq 1$};\\
         x & \mbox{if $ -1 < x < 1$};\\
         -1 & \mbox{if $x \leq -1$}.\end{array} \right. \] 
We will first solve the problem on a real line. Let $v(x)$ be a bounded and continuous function. We will solve the differential equation
\begin{equation}
\phi(x) \frac{\partial u(x)}{\partial x} = v(x)
\end{equation}
where $u$ is a bounded function.

Let $v(x)=1$. When $-1<x<1$,  we have
\begin{align*}
&x\frac{\partial u(x)}{\partial x} =v(x)\\
& u(x)=\int \frac{1}{x}dx= \ln x
\end{align*}
Obviously, $u(x)= \ln x$ is not bounded around the origin. So this differential equation is not solvable.

Then we multiply $x$ to both sides of equation (59) and get
\begin{equation}
\phi(x) \frac{\partial xu(x)}{\partial x} = xv(x)
\end{equation}
Let $v(x)=1$. When $-1<x<1$,  we have
\begin{align*}
&x\frac{\partial xu(x)}{\partial x} =xv(x),\\
& xu(x)=\int 1dx=  x,\\
&u(x)=1.
\end{align*}
When $x \geq 1$, we have
\begin{align*}
&\frac{\partial xu(x)}{\partial x} =xv(x),\\
& xu(x)=\int xdx=  \frac{1}{2} x^2,\\
& u(x)=\frac{x}{2}.
\end{align*}
Since $u(x)=\frac{x}{2}$ is not bounded, this differential equation is not solvable either.

Then we multiply $e^x$ to both sides of equation (60) and get
\begin{equation}
\phi(x) \frac{\partial xe^xu(x)}{\partial x} = xe^xv(x)
\end{equation}
Let $v(x)=1$. When $-1<x<1$,  we have
\begin{align*}
&x\frac{\partial xe^xu(x)}{\partial x} =xe^xv(x),\\
& xe^xu(x)=\int_0^x e^tdt= e^x-1,\\
&u(x)=\frac{e^x-1}{xe^x}.
\end{align*}
Since
\begin{equation*}
\lim_{x \rightarrow 0} \frac{e^x-1}{x} =1,
\end{equation*}
then when $x \rightarrow 0$, $u(x) = 1$. So $u(x)$ is bounded between -1 and 1. 

When $x \geq 1$, we have
\begin{align*}
&\frac{\partial xe^xu(x)}{\partial x} =xe^xv(x),\\
& xe^xu(x)=\int_0^x te^tdt=(te^t)|_0^x-\int_0^x e^tdt=xe^x-e^x+1,\\ 
& u(x)=\frac{xe^x-e^x+1}{ xe^x} \leq \frac{3xe^x}{ xe^x}=3.
\end{align*}
Since $u(x) \leq 3$, we know $u(x)$ is bounded when $x \geq 1$.

When $x \leq -1$, we have
\begin{align*}
&\frac{\partial xe^xu(x)}{\partial x} =xe^xv(x),\\
& xe^xu(x)=\int_0^x te^tdt=(te^t)|_0^x-\int_0^x e^tdt=xe^x-e^x+1,\\ 
& u(x)=\frac{xe^x-e^x+1}{ xe^x} \leq \frac{2xe^x}{ xe^x} + \frac{1}{xe^x}=2+\frac{1}{xe^x}.
\end{align*}
Since $\frac{1}{xe^x}$ is not bounded when $x \rightarrow -\infty$, $u(x)$ is not bounded when $x \leq 1$. 

Then we modify equation (61) into
\begin{equation}
\phi(x) \frac{\partial xe^{|x|}u(x)}{\partial x} = xe^{|x|}v(x)
\end{equation}
When $x > 0$, we have the same bounded function $u(x)$ as what we got from equation (60). When $-1 < x < 0$, we have
\begin{align*}
&x\frac{\partial xe^{-x}u(x)}{\partial x} =xe^{-x}v(x),\\
& xe^{-x}u(x)=\int_0^x e^{-t}dt= -e^{-x}+1,\\
&u(x)=\frac{-e^x+1}{xe^{-x}}.
\end{align*}
Since
\begin{equation*}
\lim_{x \rightarrow 0} \frac{-e^x+1}{x} =1,
\end{equation*}
then when $x \rightarrow 0$, $u(x) = 1$. So $u(x)$ is bounded between -1 and 0. 
When $x \leq -1$,  we have
\begin{align*}
&\frac{\partial xe^{-x}u(x)}{\partial x} =xe^{-x}v(x),\\
& xe^{-x}u(x)=\int_0^x te^{-t}dt=(-te^{-t})|_0^x-\int_0^x -e^{-t}dt=-xe^{-x}-e^{-x}+1,\\ 
& u(x)=\frac{-xe^{-x}-e^{-x}+1}{ xe^{-x}} \leq \frac{3xe^x}{ xe^x}=3
\end{align*}
Since $u(x) \leq 3$, we know $u(x)$ is bounded when $x \leq -1$. Therefore, $u(x)$ is bounded through the real line and differential equation (62) is solvable.

Now we will derive the general formulae for $u(x)$ in equation (62). 

When $0 \leq x \leq 1$, we have
\begin{align*}
&x\frac{\partial xe^xu_1(x)}{\partial x} =xe^xv(x),\\
& xe^xu_1(x)=\int_0^x e^tv(t)dt,\\
&u_1(x)=\frac{1}{x}e^{-x}\int_0^x e^tv(t)dt.
\end{align*}
When $x \geq 1$, we have
\begin{align*}
&\frac{\partial xe^xu_2(x)}{\partial x} =xe^xv(x),\\
& xe^xu_2(x)=\int_1^x te^tv(t)dt,\\ 
& u_2(x)=\frac{1}{x}e^{-x}\int_1^x [te^tv(t)dt+eu_1(1)]
\end{align*}
When $-1 \leq x \leq 0$, we have
\begin{align*}
&x\frac{\partial xe^{-x}u_3(x)}{\partial x} =xe^{-x}v(x),\\
& xe^{-x}u_3(x)=\int_0^x e^{-t}v(t)dt,\\
&u_3(x)=\frac{1}{x}e^{x}\int_0^x e^{-t}v(t)dt.
\end{align*}
When $x \leq -1$, we have
\begin{align*}
&\frac{\partial xe^{-x}u_4(x)}{\partial x} =xe^{-x}v(x),\\
& xe^{-x}u_4(x)=\int_{-1}^x te^{-t}v(t)dt,\\ 
& u_4(x)=\frac{1}{ x}e^{x} \int_{-1}^x [te^{-t}v(t)dt-eu_3(-1)]
\end{align*}

Then we are going to show the continuity of $u(x)$ if $v(x)$ is continuous by show that $u(x)$ is continuous in each of the four parts and between neighboring  parts.

Let $\epsilon$ be an arbitrary small number. 

When $0 \leq x \leq 1$, we will show that there exist a real number $\delta_1$ such that $\abs{xe^xu_1(x)-xe^xu_1(x')} < \epsilon$ holds for every pair of $(x, x')$ if $\abs{x-x'} < \delta_1$.
\begin{eqnarray*}
&&\abs{xe^xu_1(x)-x'e^{x'}u_1(x')}\\
&=& \abs{\int_0^x e^tv(t)dt-\int_0^{x'} e^tv(t)dt}\\
&=&\abs{\int_{x'}^x e^tv(t)dt}
\end{eqnarray*}
Since $v(x)$ is bounded and $e^tv(t) \leq ev(t)$ between 0 and 1, we can find a maximum value $M$ for $e^tv(t)$. Therefore, there exists $\delta_1$ such that
\begin{equation*}
\abs{\int_{x'}^x e^tv(t)dt} \leq \abs{\int_{x'}^x Mdt} < \epsilon
\end{equation*}
if $\abs{x-x'} < \delta_1$.

When $x \geq 1$, we will show that there exist a real number $\delta_2$ such that $\abs{xe^xu_2(x)-xe^xu_2(x')} < \epsilon$ holds for every pair of $(x, x')$ if $\abs{x-x'} < \delta_2$.
\begin{eqnarray*}
&&\abs{xe^xu_2(x)-x'e^{x'}u_2(x')}\\
&=& \abs{\int_1^x te^tv(t)dt-\int_1^{x'} te^tv(t)dt}\\
&=&\abs{\int_{x'}^x te^tv(t)dt}
\end{eqnarray*}
We fix $x$ at a random point, we have a finite value of $te^t$. Therefore, we have
\begin{equation*}
\abs{xe^xu_2(x)-x'e^{x'}u_2(x')} < \epsilon
\end{equation*} 
if $\delta_2$ is small enough and $\abs{x-x'} < \delta_2$.

When $-1 \leq x \leq 0$, we will show that there exist a real number $\delta_3$ such that $\abs{xe^{-x}u_3(x)-xe^{-x}u_3(x')} < \epsilon$ holds for every pair of $(x, x')$ if $\abs{x-x'} < \delta_3$.
\begin{eqnarray*}
&&\abs{xe^{-x}u_3(x)-x'e^{-x'}u_3(x')}\\
&=& \abs{\int_0^x e^{-t}v(t)dt-\int_0^{x'} e^{-t}v(t)dt}\\
&=&\abs{\int_{x'}^x e^{-t}v(t)dt}
\end{eqnarray*}
Since $v(x)$ is bounded and $e^{-t}v(t) \leq ev(t)$ between 0 and 1, we can find a maximum value $M$ for $e^{-t}v(t)$. Therefore, there exists $\delta_3$ such that
\begin{equation*}
\abs{\int_{x'}^x e^{-t}v(t)dt} \leq \abs{\int_{x'}^x Mdt} < \epsilon
\end{equation*}
if $\abs{x-x'} < \delta_3$.

When $x \leq -1$, we will show that there exist a real number $\delta_4$ such that $\abs{xe^{-x}u_4(x)-xe^{-x}u_4(x')} < \epsilon$ holds for every pair of $(x, x')$ if $\abs{x-x'} < \delta_4$.
\begin{eqnarray*}
&&\abs{xe^{-x}u_4(x)-x'e^{-x'}u_4(x')}\\
&=& \abs{\int_{-1}^x te^{-t}v(t)dt-\int_{-1}^{x'} te^{-t}v(t)dt}\\
&=&\abs{\int_{x'}^x te^{-t}v(t)dt}
\end{eqnarray*}
We fix $x$ at a random point, we have a finite value of $te^{-t}$. Therefore, we have
\begin{equation*}
\abs{xe^xu_4(x)-x'e^{x'}u_4(x')} < \epsilon
\end{equation*} 
if $\delta_4$ is small enough and $\abs{x-x'} < \delta_4$.

To show the continuity between different parts of $u(x)$, we need to have the following equations:
\begin{align}
&u_1(0)=u_3(0)\\
&u_1(1)=u_2(1)\\
&u_3(-1)=u_4(-1)
\end{align}
\begin{align*}
&\lim_{x \rightarrow 0}u_1(0)=0\\
&\lim_{x \rightarrow 0}u_3(0)=0\\
&u_2(1)=\frac{1}{e} \int_1^1 te^tv(t)dt +\frac{1}{e}eu_1(1)=u_1(1)\\
&u_4(-1)=-\frac{1}{e} \int_{-1}^{-1} te^{-t}v(t)dt +\frac{1}{e}eu_3(-1)=u_3(-1)\\
\end{align*}
Since we have shown $u(x)$ is continuous on each part and equations (63) (64) and (65) hold, we know that $u(x)$ is continuous.

\subsection{The circle case}

Consider a vector field $X=\phi(\theta)\frac{\partial}{\partial \theta}$ on circle with singularities on it. Notice that there are infinitely many singularities on the circle since we encounter one or more singularities every $2\pi$ distance. So the function $\phi(\theta)$ has to be a periodic function with period of $2\pi$. let 
\begin{equation*}
\phi(\theta)=\sin \theta
\end{equation*}
Based on previous experience, we have our differential equation in the form of
\begin{equation*}
\sin \theta \frac{d}{d\theta}[f(\theta)u(\theta)]=f(\theta)v(\theta)
\end{equation*}
where $v(\theta)$ is a continuous and periodic function, $f(\theta)$ is continuous and $u(\theta)$ is continuous, periodic and differentiable along $\theta$. Obviously, there are singularities at $\theta=0$ and $\theta=\pi$.
Now we generate the general formulae for $u(\theta)$ which is
\begin{equation*}
u(\theta)=\frac{1}{f(\theta)} \int_{-\infty}^\theta \frac{f(x)v(x)}{\sin x}dx
\end{equation*}
Since $v(\theta)$ periodic, let $x_1=x-2\pi$, we have
\begin{eqnarray*}
u(\theta+2\pi)&=&\frac{1}{f(\theta+2\pi)} \int_{-\infty}^{\theta+2\pi} \frac{f(x)v(x)}{\sin(x)}dx\\
&=&\frac{1}{f(\theta+2\pi)} \int_{-\infty}^\theta  \frac{f(x_1+2\pi)v(x_1+2\pi)}{\sin( x_1)}dx\\
&=&\frac{1}{f(\theta+2\pi)} \int_{-\infty}^\theta  \frac{f(x_1+2\pi)v(x_1)}{\sin( x_1)}dx\\
\end{eqnarray*}
So to make $u(\theta)$ to be periodic, $f(\theta)$ has to satisfy
\begin{equation*}
f(\theta+2\pi)=Cf(\theta)
\end{equation*}
where $C$ is a constant.

We found that if the function $f(x)$ cancels the singularity of $\sin x$ at 0, i.e. $\lim_{x \rightarrow 0} = 0$, then there would be a problem at $\theta =0$ since there is a $\frac{1}{f(\theta)}$ outside the integral. 

Then, to avoid the point $x=0$, we change our integral interval to $(0, \theta)$ and get
\begin{eqnarray*}
u(\theta+2\pi)&=&\frac{1}{f(\theta+2\pi)} \int_0^{\theta+2\pi} \frac{f(x)v(x)}{\sin(x)}dx\\
&=&\frac{1}{f(\theta+2\pi)} \int_0^\theta  \frac{f(x_1+2\pi)v(x_1+2\pi)}{\sin( x_1)}dx\\
&=&\frac{1}{f(\theta+2\pi)} \int_{-2\pi}^\theta  \frac{f(x_1+2\pi)v(x_1)}{\sin( x_1)}dx
\end{eqnarray*}

Since we can always find a $v(x)$ such that
\begin{equation*}
\int_{-2\pi}^0 \frac{f(x)v(x)}{\sin x} \neq 0,
\end{equation*}
there is no way for $u(\theta)$ to be periodic.

\section{Line Bundles}

\subsection{The Circle Case.}
We partition a circle into upper part and lower part with two overlaps at the left and right ends respectively. Let $S_u$ be the set that contains the upper part of the circle and $S_l$ contains the lower part. We define $L$ as a set of contains left overlap area and $R$ as the set contains the right overlap area. Both $L$ and $R$ are intersection of sets $S_u$ and $S_l$. Base on the assumptions, we have
\begin{align*}
&S_u \times \mathbb{R}(t_u)=\{ (a, t_u)\}, \mbox{ a is a point on the upper circle and } a \in S_u, t_u \in \mathbb{R}\\
&S_l \times \mathbb{R}(t_l)=\{ (a, t_l)\}, \mbox{ a is a point on the lower circle and } a \in S_l, t_l \in \mathbb{R}
\end{align*}
\begin{align*}
&L \times \mathbb{R} = \{ (a, t_u)| a=b, t_u=t_l \}\\
&R \times \mathbb{R} = \{ (a, t_u)| a=b, t_u=\frac{1}{e^{2\pi}}t_l \}
\end{align*}
Then we separate the differential equation 
\begin{equation*}
\frac{du}{d\theta}=v
\end{equation*}
into the following two equations:
\begin{equation*}
\frac{du_u}{d\theta}=v_u \mbox{,   } \frac{du_l}{d\theta}=v_l
\end{equation*}
Thus we separate functions $u$ and $v$ into four functions: $\{ u_u, u_l\}$ and $\{ v_u, v_l\}$.

On the set $S_l$ which is the lower part of the circle, we have
\begin{equation*}
t_l=u_l \mbox{ in } \{ u_u, u_l\} \mbox{ and } t_l=v_l \mbox{ in } \{ v_u, v_l\}
\end{equation*}
On the set $S_u$ which is the upper part of the circle, we have
\begin{equation*}
t_u=u_u \mbox{ in } \{ u_u, u_l\} \mbox{ and } t_u=v_u \mbox{ in } \{ v_u, v_l\}
\end{equation*}
On the set $L$ which is the left overlap part of the circle, we have
\begin{equation*}
u_u=u_l \mbox{,   } v_u=v_l 
\end{equation*}
On the set $R$ which is the right overlap part of the circle, we have
\begin{equation*}
u_u=\frac{1}{e^{2\pi}} u_l \mbox{,   } v_u=\frac{1}{e^{2\pi}} v_l 
\end{equation*}
Let
\begin{align*}
&\tilde{u_u}=\frac{u_u}{e^\theta_1} \mbox{,  } \tilde{v_u}=\frac{v_u}{e^\theta_1} \mbox{ where } -\epsilon \leq \theta_1 \leq \epsilon + \pi\\
&\tilde{u_l}=\frac{u_l}{e^\theta_2} \mbox{,  } \tilde{v_l}=\frac{v_l}{e^\theta_2} \mbox{ where } -\epsilon+ \pi \leq \theta_2 \leq \epsilon + 2\pi
\end{align*}
Notice that in set $L$, $\theta_2=\theta_1$ and in set $R$, $\theta_2=\theta_1 +2\pi$.
We have to make sure that $\tilde{u_l}$ and $\tilde{u_u},  \tilde{v_l}$ and $\tilde{v_u}$ match each other on the overlap part.

On set $L$, we have $\tilde{u_l} = \tilde{u_u}$.

On set $R$, we have
\begin{equation*}
\tilde{u_u} = \frac{u_u}{e^\theta_1} = \frac{u_l}{e^{\theta_1+2\pi}}=\frac{u_l}{e^{\theta_2}}=\tilde{u_l}.
\end{equation*}
Similarly, we have $\tilde{v_l} = \tilde{v_u}$ on both of the overlap area.

Thus our differential equation becomes
\begin{equation}
\frac{d(e^\theta \tilde{u}(\theta))}{d\theta}=e^\theta \tilde{v}(\theta)
\end{equation}
Based on previous experience, $\tilde{u}(\theta)$ is solvable and the general formulae for $\tilde{u}(\theta)$ is 
\begin{equation*}
\tilde{u}(\theta)=e^{-\theta} \int_{-\infty}^{\theta} e^t\tilde{v}(t)dt
\end{equation*}
Since functions $\tilde{u}(\theta)$ and  $\tilde{v}(\theta)$ are defined on the circle, they have to be both periodic with period of $2\pi$. So we have $ \tilde{v}(\theta)$ is periodic and we will show that $\tilde{u}(\theta)$ is periodic too. Let $t'=t-2\pi$, we have
\begin{eqnarray*}
\tilde{u}(\theta+2\pi)&=&e^{-\theta-2\pi} \int_{-\infty}^{\theta+2\pi} e^t\tilde{v}(t)dt\\
&=&e^{-\theta-2\pi} \int_{-\infty}^{\theta} e^{t'+2\pi}\tilde{v}(t'+2\pi)dt'\\
&=&e^{-\theta} \int_{-\infty}^{\theta} e^{t'}\tilde{v}(t')dt'\\
&=&\tilde{u}(\theta)
\end{eqnarray*}
Hence we know that $\tilde{u}(\theta)$ is a periodic function.

Then we show that $\tilde{u}(\theta)$ is a continuous function if $tilde{v}(\theta)$ is continuous.

Let $\epsilon$ be an arbitrary small number. We will show that there exist a real number $\delta$ such that $\abs{e^\theta\tilde{u}(\theta)-e^\theta\tilde{u}(\theta')}<\epsilon$ holds for every pair of $(\theta, \theta')$ if $\abs{\theta-\theta'} < \delta$.
\begin{eqnarray*}
&&\abs{e^\theta\tilde{u}(\theta)-e^\theta\tilde{u}(\theta')}\\
&=&|\int_{-\infty}^{\theta} e^t\tilde{v}(t)dt-\int_{-\infty}^{\theta'} e^t\tilde{v}(t)dt|\\
&=& |\int_\theta^{\theta'}e^t\tilde{v}(t)dt|
\end{eqnarray*}
We fix $\theta$ at a random point, we have a finite value of $e^t\tilde{v}(t)$. Therefore, we have
\begin{equation*}
\abs{e^\theta\tilde{u}(\theta)-e^\theta\tilde{u}(\theta')}< \epsilon
\end{equation*} 
if $\delta$ is small enough and $\abs{\theta-\theta'} < \delta$.

\subsection{Annulus with Spirals Case.}
In this section, we have two concentric circles with radius 1 and 2 respectively, and infinite many spirals between them just as what we have in section 4. The spirals goes though every point between the two circles and each of them goes to approximate both of the circles. The spirals  will neither touch each other nor the two circles so that if we pick an arbitrary small part of the foliation, we can stretch the curves to be straight parallel lines. A continuous and bounded function $V(x, y)$ is defined on the two circles and the spiral, i.e. it is defined on the closed area between the two circles. We want to find another continuous and bounded  function $U(x, y)$ which is differentiable along the spirals such that $V(x, y)=XU(x, y)$ where $X$ is the vector field. $X$ has the form of
\begin{equation*}
X = F(x, y) \frac{\partial}{\partial x}+ G(x, y) \frac{\partial}{\partial y},
\end{equation*}
where $F(x, y)$ and $G(x, y)$ are continuous functions of $x$ and $y$.

Since the definition of the functions are on circles and spirals, it would be more convenient if we transform the coordinates into polar systems. 
Let $r$, $\theta$ be two variables in the corresponding polar system. Then we have
\begin{align}
 x & =r \cos \theta \\
 y & =r \sin \theta
\end{align}
Notice that the radius $r$ is a function of the angle $\theta$, we have the expression of $r$ on spiral s:
\begin{equation}
r(\theta, s)=\frac{3}{2}+\arctan (\theta+s), 
\end{equation}
where $ -\pi \leq s \leq  \pi$. 

According to the results in section 4, we have the following transformations:

From Cartesian to Polar coordinates:
\begin{align*}
& V(x, y)=V(r \cos \theta, r \sin \theta) = v(r, \theta) \\
& U(x, y)=U(r \cos \theta, r \sin \theta) = u(r, \theta) 
\end{align*} 
From Polar to Cartesian coordinates:
\begin{align*}
& v(r, \theta)=v(\sqrt{x^2+y^2}, \arctan \frac{y}{x}) = V(x, y) \\
& u(r, \theta)=u(\sqrt{x^2+y^2}, \arctan \frac{y}{x}) = U(x, y)
\end{align*} 
The transformation of vector field $X$:
\begin{align*}
&X = F(x, y) \frac{\partial}{\partial x}+ G(x, y) \frac{\partial}{\partial y},\\
&F(x, y) = \cos \theta \frac{dr}{d\theta} -r \sin \theta,\\
&G(x, y)= \sin \theta \frac{dr}{d\theta} +r \cos \theta,\\
&V(x, y) = (F(x, y)\frac{\partial}{\partial x}+G(x, y)\frac{\partial}{\partial y})U(x, y).
\end{align*}

On the circles, we have
\begin{equation*}
v(1, \theta)=v(1, \theta+2\pi),
\end{equation*}
\begin{equation*}
v(2, \theta)=v(2, \theta+2\pi),
\end{equation*}
\begin{equation*}
u(1, \theta) = u(1, \theta +2\pi), 
\end{equation*}
\begin{equation*}
u(2, \theta) = u(2, \theta +2\pi).
\end{equation*}
On spiral s, function $v(r(\theta, s), \theta)$ has following property on the spiral:
\begin{equation}
\lim_{\theta \rightarrow -\infty} \abs{ v(r(\theta, s), \theta)-v(1, \theta)} = 0,
\end{equation}
\begin{equation}
\lim_{\theta \rightarrow +\infty} \abs{ v(r(\theta, s), \theta)-v(2, \theta)} = 0.
\end{equation}
We want to show that function $u(r(\theta, s), \theta)$ has the same property:
\begin{equation*}
\lim_{\theta \rightarrow -\infty} \abs{ u(r(\theta, s), \theta)-u(1, \theta)} = 0,
\end{equation*}
\begin{equation*}
\lim_{\theta \rightarrow +\infty} \abs{ u(r(\theta, s), \theta)-u(2, \theta)} = 0.
\end{equation*}
We partition the annulus into upper part and lower part with two overlaps at the left and right ends respectively. Let $S_u$ be the set that contains the upper part of the annulus and $S_l$ contains the lower part. We define $L$ as a set contains left overlap area and $R$ as the set contains the right overlap area. Both $L$ and $R$ are intersection of sets $S_u$ and $S_l$. Base on the assumptions, we have
\begin{align*}
&S_u \times \mathbb{R}(t_u)=\{ (a, t_u)\}, \mbox{ a is a point on the upper annulus and } a \in S_u, t_u \in \mathbb{R}\\
&S_l \times \mathbb{R}(t_l)=\{ (a, t_l)\}, \mbox{ a is a point on the lower annulus and } a \in S_l, t_l \in \mathbb{R}
\end{align*}
\begin{align*}
&L \times \mathbb{R} = \{ (a, t_u)| a=b, t_u=t_l \}\\
&R \times \mathbb{R} = \{ (a, t_u)| a=b, t_u=\frac{1}{e^{2\pi}}t_l \}
\end{align*}
Then we separate the differential equation 
\begin{equation*}
\frac{du}{d\theta}=v
\end{equation*}
into the following two equations:
\begin{equation*}
\frac{du_u}{d\theta}=v_u \mbox{,   } \frac{du_l}{d\theta}=v_l
\end{equation*}
Thus we separate functions $u$ and $v$ into four functions: $\{ u_u, u_l\}$ and $\{ v_u, v_l\}$.

On the set $S_l$ which is the lower part of the annulus, we have
\begin{equation*}
t_l=u_l \mbox{ in } \{ u_u, u_l\} \mbox{ and } t_l=v_l \mbox{ in } \{ v_u, v_l\}
\end{equation*}
On the set $S_u$ which is the upper part of the annulus, we have
\begin{equation*}
t_u=u_u \mbox{ in } \{ u_u, u_l\} \mbox{ and } t_u=v_u \mbox{ in } \{ v_u, v_l\}
\end{equation*}
On the set $L$ which is the left overlap part of the annulus, we have
\begin{equation*}
u_u=u_l \mbox{,   } v_u=v_l 
\end{equation*}
On the set $R$ which is the right overlap part of the annulus, we have
\begin{equation*}
u_u=\frac{1}{e^{2\pi}} u_l \mbox{,   } v_u=\frac{1}{e^{2\pi}} v_l 
\end{equation*}
Let
\begin{align*}
&\tilde{u_u}=\frac{u_u}{e^\theta_1} \mbox{,  } \tilde{v_u}=\frac{v_u}{e^\theta_1} \mbox{ where } -\epsilon \leq \theta_1 \leq \epsilon + \pi\\
&\tilde{u_l}=\frac{u_l}{e^\theta_2} \mbox{,  } \tilde{v_l}=\frac{v_l}{e^\theta_2} \mbox{ where } -\epsilon+ \pi \leq \theta_2 \leq \epsilon + 2\pi
\end{align*}
Notice that in set $L$, $\theta_2=\theta_1$ and in set $R$, $\theta_2=\theta_1 +2\pi$.
We have to make sure that $\tilde{u_l}$ and $\tilde{u_u},  \tilde{v_l}$ and $\tilde{v_u}$ match each other on the overlap part.

On set $L$, we have $\tilde{u_l} = \tilde{u_u}$.

On set $R$, we have
\begin{equation*}
\tilde{u_u} = \frac{u_u}{e^\theta_1} = \frac{u_l}{e^{\theta_1+2\pi}}=\frac{u_l}{e^{\theta_2}}=\tilde{u_l}.
\end{equation*}
Similarly, we have $\tilde{v_l} = \tilde{v_u}$ on both of the overlap area.

Thus our differential equation becomes
\begin{equation}
\frac{d(e^\theta \tilde{u}(r(\theta, s), \theta)}{d\theta}=e^\theta \tilde{v}(r(\theta, s), \theta)
\end{equation}
Based on previous experience, $\tilde{u}(r(\theta, s), \theta)$ is solvable and the general formulae for $\tilde{u}(r(\theta, s), \theta)$ is 
\begin{equation*}
\tilde{u}(r(\theta, s), \theta)=e^{-\theta} \int_{-\infty}^{\theta} e^t\tilde{v}((r(t, s), t))dt.
\end{equation*}
We have shown the continuity of $\tilde{u}(r(\theta, s), \theta)$ in Section 2.4.

\subsection{General Foliations Case (Proof of Main Result).}

Consider a laminated set $K$ on $\mathbb{R}$ which is a union of integral curves. We cover the laminated set $K$ by finitely many boxes. We have an orientation of the set $K$ that define the direction that the arc length increases along the curves. Imagine there is a particle moving along the curves at a velocity of 1. Assume along the direction of the curves, there are drops in the overlap area. We define the drop between boxes as following:

In each box $B$, we have a function $f_B$. We assume that at the left boundary of the box, $f_B=0$ for all curves. That is to say, when $t=0$, we have $f_B(t)=0$ where $t$ is distance that a particle go through in the box since the velocity of the particle is always 1. 

Prove. Suppose we have two boxes $B_1$ and $B_2$ and the overlap of the two boxes $O_{12}$. There are some curves go from $B_1$ to $B_2$. On curve $L$, we have a dropping constant given by
\begin{equation*}
C_{12}^L=\frac{f_{B_1}}{f_{B_2}}
\end{equation*}
Based on the assumptions, on leaf $L$ we have
\begin{align*}
&B_1 \times \mathbb{R}(t_1) = \{(a, t_1)\}, \mbox{$a$ is a point in $B_1$ and $a \in B_1, t_1 \in \mathbb{R}$}\\
&B_2 \times \mathbb{R}(t_2) = \{(b, t_2)\}, \mbox{$a$ is a point in $B_2$ and $b \in B_2, t_1 \in \mathbb{R}$}
\end{align*}
\begin{equation*}
O_{12} \times \mathbb{R} = \{(a, t_1)|a=b, t_2=C_{12}t_1\}
\end{equation*}
Assume we have three boxes $B_1$ $B_2$ and $B_3$ and they all have overlaps with each other. Then we have
\begin{equation*}
C_{12}C_{23} = \frac{f_{B_1}}{f_{B_2}} {f_{B_2}}{f_{B_3}}={f_{B_1}}{f_{B_3}}.
\end{equation*}
This proves the consistency of our assumption. 

Then we show our main results. We separate the differential equation 
\begin{equation*}
\frac{du}{dt}=v
\end{equation*}
into (infinite) many equations:
\begin{equation*}
\frac{du_1}{dt}=v_1, \frac{du_2}{dt}=v_2, ..., \frac{du_n}{dt}=v_n.
\end{equation*}
In box $B_i$, we have
\begin{equation*}
t_i=u_i, t_i=v_i.
\end{equation*}
On overlap $O_{ij}$, we have
\begin{equation*}
u_j=u_iC_{ij}, v_j=v_iC_{ij}.
\end{equation*}
Let
\begin{equation*}
\tilde{u_i} =\frac{u_i}{e^{t_i}}, \tilde{v_i} =\frac{v_i}{e^{t_i}}, 
\end{equation*}
\begin{equation*}
\tilde{u_j} =\frac{u_j}{e^{t_j}}, \tilde{v_j} =\frac{v_j}{e^{t_j}}, 
\end{equation*}
Notice that on overlap $O_{ij}$, $t_j=t_i+\ln C_{ij}$.
\begin{equation*}
\tilde{u_j} =\frac{u_j}{e^{t_j}}=\frac{u_iC_{ij}}{e^{t_j}}=\frac{u_i}{e^{t_i}}=\tilde{u_i}
\end{equation*}
Similarly, we have $\tilde{v_i}=\tilde{v_j}$ on the overlap $O_{ij}$.

Therefore, our differential equation becomes
\begin{equation*}
\frac{d(e^t\tilde{u})}{dt}=e^t\tilde{v}.
\end{equation*}
which is the equation we start with this paper. The general formulae of $\tilde{u}$ is
\begin{equation*}
\tilde{u}(t)=e^{-t}\int_{\infty}^{t} \tilde{v}(w)dw
\end{equation*}
We have shown the continuity of $\tilde{u}$ in section 2.

\bigskip

\noindent John E. Fornaess\\
Mathematics Department\\
The University of Michigan\\
 Ann Arbor, Michigan 48109\\
USA\\
fornaess@umich.edu\\

\noindent Xiaoai Chai\\
Mathematics Department\\
The University of Michigan\\
 Ann Arbor, Michigan 48109\\
USA\\
xachai@umich.edu\\

\end{document}